\documentclass[12pt, a4paper]{amsart}
\usepackage[dvips]{epsfig}
\usepackage{amsmath}
\usepackage{amssymb}
\usepackage{amsfonts}
\usepackage{dsfont}
\usepackage{bbm}
\usepackage{amsthm}
\usepackage{amsbsy}
\usepackage{amsgen}
\usepackage{amscd}
\usepackage{amsopn}
\usepackage{amstext}
\usepackage{amsxtra}
\usepackage[utf8]{inputenc}
\usepackage{enumerate}
\usepackage{hyperref}
\usepackage{lineno}
\usepackage{marginnote}
\usepackage{verbatim}
\usepackage{calrsfs}
\usepackage{times, graphicx}
\usepackage{eso-pic}
\usepackage{lastpage}
\DeclareMathAlphabet{\pazocal}{OMS}{zplm}{m}{n}
\usepackage[foot]{amsaddr}

\numberwithin{equation}{section}

\makeatletter
\newcommand*\rel@kern[1]{\kern#1\dimexpr\macc@kerna}
\newcommand*\widebar[1]{%
  \begingroup
  \def\mathaccent##1##2{%
    \rel@kern{0.8}%
    \overline{\rel@kern{-0.8}\macc@nucleus\rel@kern{0.2}}%
    \rel@kern{-0.2}%
  }%
  \macc@depth\@ne
  \let\math@bgroup\@empty \let\math@egroup\macc@set@skewchar
  \mathsurround\z@ \frozen@everymath{\mathgroup\macc@group\relax}%
  \macc@set@skewchar\relax
  \let\mathaccentV\macc@nested@a
  \macc@nested@a\relax111{#1}%
  \endgroup
}
\makeatother

\newtheorem{theorem}{Theorem}[section]
\newtheorem{proposition}[theorem]{Proposition}
\newtheorem{lemma}[theorem]{Lemma}

\theoremstyle{definition}
\newtheorem{definition}[theorem]{Definition}

\newcommand{\vol}{\operatorname{Vol}}

\newcommand{\prob}{\pazocal{P}r}

\newcommand{\Dc}{\mathcal{D}}

\newcommand{\Sc}{\pazocal{S}}
\newcommand{\Nc}{\pazocal{N}}
\newcommand{\Ec}{\pazocal{E}}

\newcommand{\Fc}{\pazocal{F}}
\newcommand{\Kc}{\pazocal{K}}

\newcommand {\E} {\mathbb{E}}
\newcommand {\M} {\pazocal{M}}
\newcommand {\R} {\mathbb{R}}

\newcommand {\Z} {\mathbb{Z}}

\newcommand {\Cc} {\pazocal{C}}
\newcommand {\Lc} {\mathcal{L}}

\newcommand {\Pc} {\mathcal{P}}

\newcommand {\C} {\mathbb{C}}
\newcommand {\Ccc} {\mathcal{C}}

\newcommand {\Ac} {\pazocal{A}}

\newcommand {\Tb} {\mathbb{T}}

\newcommand {\Vc} {\mathcal{V}}

\newcommand{\id}{\mathds{1}}

\graphicspath{{../pics/}}

\begin{document}

\title[Mean conservation for Gaussian ensembles]{Mean conservation of nodal volume and connectivity measures for Gaussian ensembles}
\author{Dmitry Beliaev\textsuperscript{1}}
\email{dmitry.belyaev@maths.ox.ac.uk}
\address{\textsuperscript{1}Mathematical Institute, University of Oxford}
\author{Stephen Muirhead\textsuperscript{2}}
\email{s.muirhead@qmul.ac.uk}
\address{\textsuperscript{2}School of Mathematical Sciences, Queen Mary University of London}
\author{Igor Wigman\textsuperscript{3}}
\email{igor.wigman@kcl.ac.uk}
\address{\textsuperscript{3}Department of Mathematics, King's College London}
\dedicatory{Dedicated to the memory of Jean Bourgain}

\date{\today}
\begin{abstract}
We study in depth the nesting graph and volume distribution of the nodal domains of a Gaussian field, which have been shown in previous works to exhibit asymptotic laws. A striking link is established between the asymptotic mean connectivity of a nodal domain (i.e.\ the vertex degree in its nesting graph) and the positivity of the percolation probability of the field, along with a direct dependence of the average nodal volume on the percolation probability. Our results support the prevailing ansatz that the mean connectivity and volume of a nodal domain is conserved for generic random fields in dimension $d=2$ but not in $d \ge 3$, and are applied to a number of concrete motivating examples.

\end{abstract}
	
\maketitle

\section{Introduction}

\subsection{The connectivity measure for nodal domains of Euclidean Gaussian fields}

\label{sec:conn meas Euclid}
\begin{figure}[ht]
\centering
\includegraphics[width=0.9\textwidth]{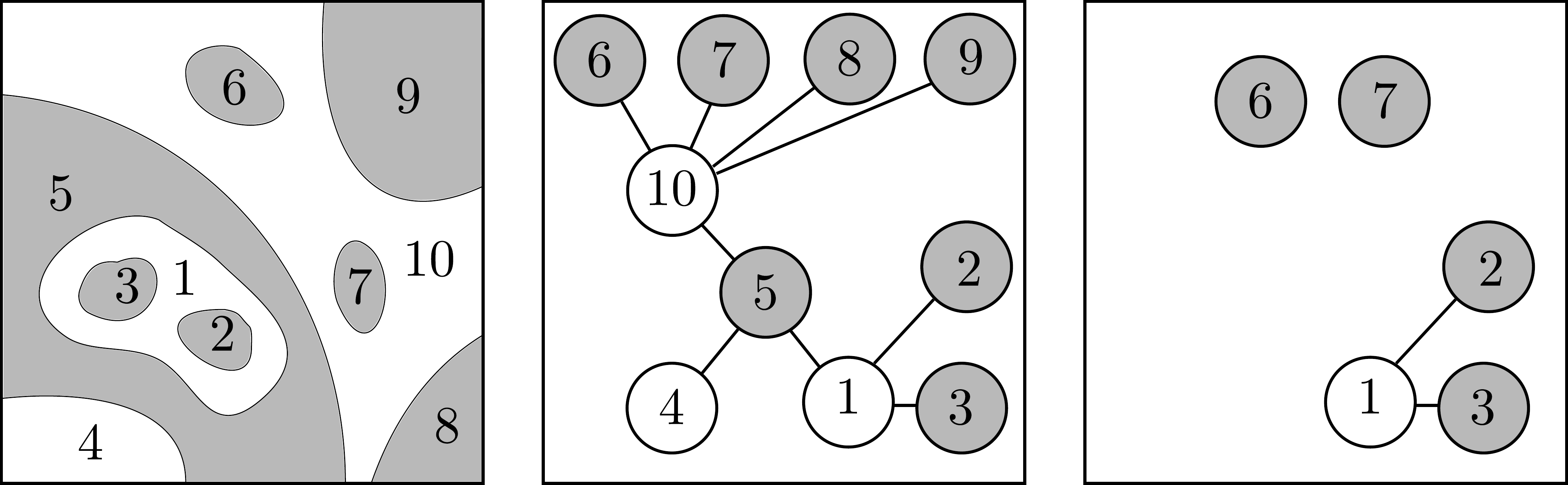}
\caption{
Left: A sketch of the nodal domains of a functions (grey and white are positive and negative domains respectively). Middle: The connectivity graph of the nodal domains; note that some of the domains might connect outside of the box, but this graph takes into account only the domains of the restricted function. Right: the corresponding graph $G(R)$, whose vertices correspond to nodal domains that are completely inside the box; typically, this graph is not a tree but a forest.}
\label{fig:nesting tree}
\end{figure}

Let $d\ge 2$ and $F:\R^{d}\rightarrow\R$ be a centred stationary $C^3$-smooth Gaussian field. We are interested in the topological structure of the nodal set $\Ac(F):=F^{-1}(0)$, of high importance in various disciplines including oceanography \cite{Longuet}, engineering \cite{Rice,Swerling} and cosmology (see for example \cite{Bond, Park} and the references therein).
A {\em nodal component} of $F$ is a connected component of $\Ac(F)$, and a {\em nodal domain} is a connected component of the complement $\R^{d}\setminus \Ac(F)$. We encode the topological structure of $\Ac(F)$ as follows: let $\Omega_{F}$ be the (a.s.\ locally finite) collection of nodal domains of~$F$, let $\Cc_{F}$ be the collection of nodal components of $F$, and define the {\em nesting graph} $G=G(F)=(V,E)$ with vertex set $V=V(F)=\Omega_{F}$ and edge set $E=E(F)=\Cc_{F}$ so that two domains $v_{1},v_{2}\in V$ are adjacent in $G$ via $e\in E$ if the corresponding domains share $e$ as a common boundary component. Figure~\ref{fig:nesting tree} exhibits a fragment of the nesting graph $G$ for some sample function $F_{\omega}$.

\smallskip
Sarnak and Wigman ~\cite{SaWi} studied the nesting graph $G$ for a generic stationary Gaussian field $F$. Observe that, by Jordan's Theorem, $G$ is a.s.\ an (infinite) tree, and so the structure of $G$ is largely encapsulated by the degrees of the vertexes $V$ (`the connectivity of the nodal domains'). Under very mild extra assumptions (to be given in \S\ref{sec:mean conn, Euclid}), Sarnak--Wigman established a law of large numbers for the connectivity in the following sense. Let $B(R)\subseteq \R^{d}$ denote the ball of radius $R$, and let $G(R)=(V(R),E(R))$ be the restriction of $G$ to $B(R)$, i.e.\ the graph induced by restricting $G$ to the vertices $V(R)\subseteq V$ that correspond to domains that are {\em fully contained within} $B(R)$; $G(R)$ might fail to be a tree but is necessarily a collection of disjoint trees (a forest). Letting
\[ d(v)=d_{R}(v)\in\Z_{\ge 0} \]
denote the degree of $v\in V(R)$ (w.r.t.\ $G(R)$), define the {\em empirical connectivity measure}
\begin{equation}
\label{eq:conn meas F;R def}
\mu_{\Gamma(F);R} = \frac{1}{|V(R)|}\sum\limits_{v\in V(R)}\delta_{d(v)}
\end{equation}
on $\Z_{\ge 0}$. Sarnak--Wigman showed ~\cite[Theorem 3.3]{SaWi} that, for a wide class of stationary Gaussian fields~$F$, as $R\rightarrow\infty$ the (random) probability measure $\mu_{\Gamma(F);R}$ tends to a {\em deterministic} probability measure $\mu_{\Gamma(F)}$ on $\Z_{\ge 0}$ (`the limit connectivity measure') that depends on the law of~$F$. More precisely, they proved that
\[ \Dc\left(\mu_{\Gamma(F);R}, \mu_{\Gamma(F)}\right)\rightarrow 0   \]
in probability as $R\rightarrow\infty$, with $\Dc(\cdot,\cdot)$ the total variation distance on probability measures on~$\Z_{\ge 0}$ (see~\eqref{eq:D(mu1,mu2) dist def} below).

\smallskip
The properties of the limit connectivity measure $\mu_{\Gamma(F)}$ are of fundamental importance, and Sarnak--Wigman raised~\cite[p.\ 13]{SaWi} the question of the {\em mean} connectivity of the limit measure~$\mu_{\Gamma(F)}$. Since~$G$ (and hence $G(R)$) contains no cycles, the mean of the empirical connectivity measures $\mu_{\Gamma(F);R}$ satisfy
\begin{equation*}
\begin{split}
\sum\limits_{k=0}^{\infty}k\cdot \mu_{\Gamma(F);R}(k) &= \frac{1}{|V(R)|}\sum\limits_{v\in V(R)} d(v) = \frac{1}{|V(R)|}\cdot 2|E(R)|\\&\le
\frac{1}{|V(R)|}\cdot 2(|V(R)|-1) = 2-\frac{2}{|V(R)|} \le 2.
\end{split}
\end{equation*}
One can deduce via Fatou's lemma ~\cite[p.\ 31]{SaWi} that
\begin{equation*}
\sum\limits_{k=0}^{\infty}k\cdot \mu_{\Gamma(F)}(k)  \le 2;
\end{equation*}
i.e.\ the mean connectivity of the limit measure $\mu_{\Gamma(F)}$ is bounded by $2$. It is then crucial to determine whether the equality
\begin{equation}
\label{eq:mean lim conn=2}
\sum\limits_{k=0}^{\infty}k\cdot \mu_{\Gamma(F)}(k) = 2
\end{equation}
holds, for if it does {\em not}, then this indicates a non-local `escape of topology' when passing to the limit.

\smallskip
Numerical experiments of Barnett--Jin (presented within ~\cite{SaWi}) seem to indicate that \eqref{eq:mean lim conn=2} {\em fails} for the monochromatic random wave and some band-limited Gaussian fields on $\R^{2}$, the motivational examples of \cite{SaWi} (see \S\ref{ssec:ens} below for definitions). On the other hand, this may be a numerical artefact due to the slow convergence of the series on the l.h.s.\ of \eqref{eq:mean lim conn=2}, reflecting the slow conjectured decay of $\mu_{\Gamma(F)}(k)$. Indeed, borrowed from percolation theory (as inspired by~\cite{BS}), it is plausible ~\cite{SaWi,KlZi} that $\mu_{\Gamma(F)}(k)$ decays only as
\begin{equation}
\label{eq:mu(k) approx 1/k^alpha}
\mu_{\Gamma(F)}(k) \approx \frac{1}{k^{\alpha}} \ , \quad  \alpha = \frac{187}{91} =2.0549\ldots,
\end{equation}
where $\alpha$ is the `Fisher exponent' that describes the {\em area} distribution of percolation clusters in Bernoulli percolation on $\Z^2$ (suggesting ~\cite{KlZi} that the connectivity of a typical domain is proportional to its area); the numerical investigations of Barnett--Jin for band-limited Gaussian fields showed consistency with \eqref{eq:mu(k) approx 1/k^alpha}, although the results were not conclusive.

\smallskip
In this manuscript we address the question of whether \eqref{eq:mean lim conn=2} holds for a wide class of smooth Gaussian fields. We believe that, contrary to Barnett--Jin's numerics, our results serve as striking evidence that~\eqref{eq:mean lim conn=2} {\em does} hold for generic fields on $\R^{2}$, including all the examples considered by Sarnak--Wigman. More precisely, our main result (Theorem \ref{thm:sum(k mu(k))=2<=>no perc} below) shows that \eqref{eq:mean lim conn=2} is essentially equivalent to the nodal domains of $F$ {\em failing to percolate}, in a sense to be made rigorous. Since, in light of~\cite{Al,BeGa,BMW,BS}, the nodal domains of a generic $F$ do not percolate if $d=2$, and, in line with ~\cite{BrLeMa,DrPrRo,Sz1}, do percolate in higher dimensions (see the discussion in \S\ref{ssec:perc} below), we believe that the equality~\eqref{eq:mean lim conn=2} holds for generic random fields on $\R^{d}$ {\em if and only if} $d=2$. To support our statement we establish this claim rigorously for a particular class of Gaussian fields on $\R^2$,  including the important special case of the Bargmann--Fock field (see \S\ref{sec:Kostlan BF ex} for details).

\subsection{Percolation probabilities for random fields}
\label{ssec:perc}
The study of the percolation of excursion sets of random fields was initiated by Molchanov--Stepanov ~\cite{MoSt1}. For a random field $$F:\R^{d}\rightarrow\R$$ and a number $u\in (-\infty,+\infty)$ one is interested whether the excursion set\footnote{In the original treatment $F^{-1}(-\infty,u)$ is studied.} $F^{-1}(u,\infty)$ {\em percolates}, i.e.\ contains an {\em unbounded component}. They found ~\cite{MoSt1} that, much like in lattice percolation, there exists a {\em critical level} $u^{*}=u^{*}(F)\in [-\infty,+\infty]$, finite or infinite, so that for $u>u^{*}$, $F^{-1}(u,+\infty)$ does not percolate a.s., whereas for $u<u^{*}$, $F^{-1}(u,+\infty)$ does percolate a.s.\ (with little information at $u=u^{*}$, although it is expected that there is no percolation at the critical level). Various criteria for when the critical level $u^{*}$ is finite were also addressed~\cite{MoSt1,MoSt2}, along with other related questions.

\smallskip
In our case one is only interested in the {\em nodal set}, being the boundary of $F^{-1}(0,+\infty)$, and so the question whether $u^{*}>0$ or $u^{*}\le 0$ is crucial; indeed if $u^* > 0$ then a.s.\ there exist giant percolating nodal domains (likely unique up to sign) that cover a positive proportion of the entire space (see, e.g.,~\cite[Theorem $16$ on p.\ 76]{BoRi}). Although this question has been resolved rigorously in only a few special cases, the picture that has emerged from the physics literature (see, e.g.,~\cite{BS}) is that $u^* = 0$ for generic centred random fields on $\R^2$, with no percolation of the nodal domains. Early work of Alexander ~\cite{Al} proved that the level lines $\{F(x)=u\}$ of a stationary-ergodic planar positive-correlated Gaussian field are a.s.\ bounded, which by the symmetry of centred Gaussian fields implies immediately that $u^\ast \le 0$. Moreover, Bogomolny--Schmidt ~\cite{BS} gave a heuristic argument demonstrating that $u^\ast = 0$ for the monochromatic random wave on the plane, essentially by comparing the random wave model to critical Bernoulli percolation on the square lattice $\Z^2$ (i.e.\ where every edge is included independently with probability $p=1/2$). Very recent results \cite{BeGa, RiVa, MuVa} have confirmed that $u^\ast=0$ for a family of planar Gaussian fields with positive and rapidly-decaying correlations, and also verified the absence of percolation of the nodal domains; an important example to which these results apply is the Bargmann--Fock field (see~\S\ref{sec:Kostlan BF ex} below).

\smallskip
On the other hand, numerical experiments recently conducted by Barnett--Jin (presented within ~\cite{SaWi}) indicate that, somewhat surprisingly, a generic centred random field $F$ on $\R^d$, $d \ge 3$, {\em does} possess giant nodal component consuming a huge proportion of the space. Sarnak ~\cite{Sa} observed that this distinction could be attributed to the fact that, for $d \ge 3$, the critical level $u^\ast$ is likely to be strictly positive, i.e.\ the nodal domains correspond to the supercritical regime for $d\ge 3$ (whereas for $d=2$ they correspond to the critical regime). This is consistent with recent results of Drewitz--Pr\'{e}vost--Rodriguez~\cite{DrPrRo} who proved that $u^\ast > 0$ for a family of {\em strongly correlated} Gaussian fields on $\Z^{3}$, including the important case of the massless harmonic crystal system with correlations decaying as
\[r_{F}(x,y)\approx \frac{1}{|x-y|}, \quad x,y\in\Z^{3}, \]
(the fact that $0 \le u^\ast < \infty$ had been previously established in \cite{BrLeMa}).

\smallskip
To state our results we make use of the following notion of {\em percolation probability}:

\begin{definition}[Percolation probability associated to a Gaussian field]
\label{def:perc prob Euclid}
Let $d\ge 2$ and $F:\R^{d}\rightarrow\R$ a $C^3$-smooth stationary Gaussian field.

\begin{enumerate}
\item For two closed sets $A,B\subseteq\R^{d}$, we define the event
\begin{equation*}
\left\{A \stackrel{F}{\longleftrightarrow} B \right\}
\end{equation*}
that there exists a nodal domain of $F$ whose closure intersects both $A$ and $B$ (`$A$ and $B$ are connected by a nodal domain of $F$'). If $A$ is the boundary of a nodal domain, then this means that there is a nodal domain \textit{adjacent} to $A$ whose closure intersects $B$.
\item The percolation probability associated to $F$ is the probability
\begin{equation*}
\Pc^{F} := \prob \left(\left\{0 \stackrel{F}{\longleftrightarrow} \infty \right\}\right)
\end{equation*}
of the event that the origin is contained in an unbounded nodal domain of $F$ (note the slight abuse of notation, where we replace a point with the corresponding singleton). Equivalently,
\begin{equation*}
\Pc^{F} := \lim\limits_{R\rightarrow\infty}\prob\left(\left\{ 0\stackrel{F}{\longleftrightarrow} \partial [-R,R]^{d}  \right\} \right)
\end{equation*}
is the limit probability of the event that the origin is contained in a nodal domain of $F$ intersecting
the boundary of a large cube $[-R,R]^{d}\subseteq\R^{d}$.
\item We say that $F$ {\em percolates} if the associated percolation probability $\Pc = \Pc^F$ is strictly positive.
\end{enumerate}
\end{definition}

It is evident that, for a continuous stationary random field, $\Pc=0$ if $u^*<0$, and $\Pc>0$ if $u^*>0$; it is moreover strongly believed, and in some cases rigorously known, that $\Pc=0$ in the case $u^*=0$. In light of the above discussion, it is natural to expect that, for a generic centred random field on $\R^d$, $\Pc > 0$ if and only if $d \ge 3$.

\subsection{Mean (non-)conservation for connectivity measures, Euclidean case}
\label{sec:mean conn, Euclid}


\smallskip
A centred continuous Gaussian field $F:\R^{d}\rightarrow\R$ is uniquely determined, via Kolmogorov's Theorem, by its covariance function
\begin{equation*}
r_{F}:\R^{d}\times\R^{d}\rightarrow\R^{d}, \quad r_{F}(x,y):= \E[F(x)\cdot F(y)] .
\end{equation*}
If $F$ is stationary then, with the usual abuse of notation,
\begin{equation*}
r_{F}(x,y) = r_{F}(0,x-y)=:r_{F}(x-y),
\end{equation*}
where now $r_{F}:\R^{d}\rightarrow\R$. Equivalently, $F$ is determined by the {\em spectral measure} $\rho = \rho_{F}$ of $F$, which is the Fourier transform of $r_{F}$; $\rho$ is a positive measure on $\R^{d}$ by Bochner's Theorem, and without loss of generality we may assume that $\rho$ is a {\em probability} measure (this corresponds to fixing $r_F(0)=1$).
As in Nazarov--Sodin ~\cite{So12,NaSo15} and Sarnak--Wigman ~\cite{SaWi}, we make the following basic assumptions on $\rho$:

\begin{definition}[Axioms on the spectral measure]
\label{def:axioms rho1-3}
$\,$
\begin{enumerate}
\item[$(\rho 1)$] The measure $\rho$ has no atoms.
\item[$(\rho 2)$] For some $p>6$,
\begin{equation*}
\int\limits_{\R^{d}} \|\lambda\|^{p} \, d\rho(\lambda) < \infty.
\end{equation*}
\item[$(\rho 3)$] The support of $\rho$ does not lie in a hyperplane in $\R^{d}$.
\end{enumerate}
\end{definition}
\noindent These assumptions imply, respectively, that $f$ is ergodic, has $C^3$-smooth sample paths a.s., and is non-degenerate.

\smallskip
Let $F:\R^{d}\rightarrow\R$ be a centred stationary Gaussian field whose spectral measure $\rho_{F}$ satisfies the axioms $(\rho 2)$--$(\rho 3)$. Nazarov--Sodin considered the total number $\Nc(F;R)$ of nodal domains of $F$ lying entirely within a large ball $B(R)$, and proved \cite{So12,NaSo15} that there exists a number $c_{NS}(\rho)\ge 0$
(`the Nazarov--Sodin constant of $F$') such that, as $R\rightarrow\infty$, we have
\[ \E[\Nc(F;R)]= c_{NS}(\rho)\cdot\vol B(R)+o_{R\rightarrow\infty}(R^{d}) ,\]
Under the additional assumption $(\rho 1)$, they moreover established the convergence in mean
\begin{equation}
\label{eq:NS conv mean}
\E\left[\left|\frac{\Nc(F;R)}{\vol B(R)}-c_{NS}(\rho)\right|\right] \rightarrow 0, \quad \text{as } R \to \infty,
\end{equation}
in particular implying a version of the law of large numbers
\begin{equation}
\label{eq:NS law large numb}
\prob\left(\left|\frac{\Nc(F;R)}{\vol B(R)}-c_{NS}(\rho)\right|> \epsilon\right) \rightarrow 0, \quad \text{as } R \to \infty,
\end{equation}
for every $\epsilon>0$. We will also need the assumption
\[  (\rho 4)\ \  c_{NS}(\rho) > 0, \phantom{aaaaaaaaaaaaaaaaaaaaaaaaaaaaaaaaaaaaaaaaaaaaaaaaaaaaa} \]
satisfied in most natural examples, which endows $\Nc(F;R)$ with the proper asymptotic scaling in~\eqref{eq:NS conv mean} (in fact, if $(\rho 4)$ fails then $\Nc(F;R) =0$ a.s.).

\smallskip
Recall the construction of the (random) nesting graph $G=G(F)=(V,E)$ corresponding to~$F$ in~\S\ref{sec:conn meas Euclid}, and recall also the restriction $G(R)$ to the ball $B(R)$ and the empirical connectivity measure $\mu_{\Gamma(F);R}$ defined in \eqref{eq:conn meas F;R def}. Under the assumptions $(\rho 1)$--$(\rho 4)$, Sarnak--Wigman ~\cite{SaWi} established the existence of a (deterministic) probability measure $\mu_{\Gamma(F)}$ on $\Z_{\ge 0}$, such that for every $\epsilon>0$,
\begin{equation}
\label{eq:d(mu(R),mu)>eps->0}
\lim\limits_{R\rightarrow\infty}\prob\left(\Dc(\mu_{\Gamma(F);R},\mu_{\Gamma(F)}) > \epsilon\right) = 0,
\end{equation}
where the distance function $\Dc(\cdot,\cdot)$ is defined as
\begin{equation}
\label{eq:D(mu1,mu2) dist def}
\Dc(\mu_{1},\mu_{2}):=\sup\limits_{A\subseteq \Z_{\ge 0}} |\mu_{1}(A)-\mu_{2}(A)|.
\end{equation}
Moreover, under a mild further condition on $\rho_{F}$, Sarnak--Wigman ~\cite{SaWi} showed that $\mu_{\Gamma(F)}$ charges the whole of $\Z_{\ge 0}$. Our first principal result asserts that, under the above assumptions on~$F$, the mean connectivity of the limit distribution $\mu_{\Gamma(F)}$ is equal to $2$ {\em if and only if} the percolation probability $\Pc$ is zero:

\begin{theorem}
\label{thm:sum(k mu(k))=2<=>no perc}
Let $F:\R^{d}\rightarrow\R$ be a continuous stationary Gaussian field whose spectral measure satisfies $(\rho 1)$--$(\rho 4)$. Let $\mu_{\Gamma(F)}$ be the limit connectivity measure defined in \eqref{eq:d(mu(R),mu)>eps->0}. Denote $\Pc=\Pc^{F}$ to be the percolation probability associated to $F$ as in Definition \ref{def:perc prob Euclid}. Then the equality
\begin{equation*}
\sum\limits_{k=0}^{\infty}k\cdot \mu_{\Gamma(F)}(k) = 2
\end{equation*}
holds if and only if $\Pc=0$.
\end{theorem}


\subsection{Mean (non-)conservation for volume distribution}
Other than the connectivity, one is also interested in the {\em empirical volume  distribution} of the nodal domains. Recall that $\Nc(F;R)$ denotes the total number of nodal domains of $F$ entirely lying inside a large ball $B(R)$. For $t>0$, let $\Nc(F,t;R)$ denote the number of such nodal domains entirely lying in $B(R)$ of volume $<t$. Refining the work of Nazarov--Sodin, Beliaev--Wigman ~\cite[Theorem 3.1]{BeWi} established that, under the same assumptions $(\rho 1)$--$(\rho 4)$ on $\rho$, the empirical volume distribution $\Nc(F,t;R)$ obeys a law of large numbers. More precisely, there exists a (deterministic) cumulative distribution function $\Psi_{F}:\R_{>0} \rightarrow [0,1]$ such that, for all continuity points $t$ of $\Psi_{F}(\cdot)$,
\begin{equation}
\label{eq:PsiF lim vol mean}
\E\left[\left| \frac{\Nc(F,t;R)}{c_{NS}(\rho)\cdot \vol B(R)} - \Psi_{F}(t)   \right|\right] \rightarrow 0
\end{equation}
as $R\rightarrow\infty$, or equivalently (in light of \eqref{eq:NS conv mean}),
\begin{equation*}
\E\left[\left| \frac{\Nc(F,t;R)}{\Nc(F;R)}\cdot\mathbbm{1}_{\Nc(F;R)>0} - \Psi_{F}(t) \right|\right] \rightarrow 0;
\end{equation*}
(the indicator controls the negligible event that $B(R)$ contains no nodal domains).

\smallskip
Similarly to the question as to whether \eqref{eq:mean lim conn=2} holds for the limit connectivity distribution, one is also interested in the mean volume of the limit distribution $\Psi_{F}$:
\begin{equation}
\label{eq:mean vol}
\int\limits_{0}^{\infty}\left( 1-\Psi_{F}(t) \right)dt.
\end{equation}
Bearing in mind that, in light of Nazarov--Sodin's \eqref{eq:NS conv mean}, the `empirical mean volume' should be about
\begin{equation*}
\frac{\vol B(R)}{\Nc(F;R)} \approx \frac{\vol B(R)}{c_{NS}(\rho)\cdot \vol B(R)} = \frac{1}{c_{NS}(\rho)},
\end{equation*}
one might expect the mean \eqref{eq:mean vol} to be equal to $\frac{1}{c_{NS}(\rho)}$. 
As such, one wishes to verify whether the triple equality 
\begin{equation}
\label{eq:mean vol lim Eucl=1/cNS}
\int\limits_{0}^{\infty}\left( 1-\Psi_{F}(t) \right)dt = \frac{1}{c_{NS}(\rho)} =\lim\limits_{R\rightarrow\infty}\E\left[\frac{\vol B(R)}{\Nc(F;R)}
\cdot\mathbbm{1}_{\Nc(F;R)>0}\right]
\end{equation}
holds (where the convergence of the r.h.s.\ of \eqref{eq:mean vol lim Eucl=1/cNS} is understood {\em in mean}); for if \eqref{eq:mean vol lim Eucl=1/cNS} fails, then this indicates an `escape of mass' in the limit. Our second main result (Theorem \ref{thm:mean cons vol Eucl} below) again verifies a connection between the `escape of mass' and the percolation probability~$\Pc$. Indeed, compared to Theorem \ref{thm:sum(k mu(k))=2<=>no perc}, we are able to explicitly quantify the failure of \eqref{eq:mean vol lim Eucl=1/cNS} as a function of the percolation probability (see \eqref{eq:mean vol lim=(1-perc)/cNS}). To state our result in full, we need to introduce one further assumption on the Gaussian field $F$:

\begin{definition}[Nodal lower concentration]
\label{def:nodal low conc}
Suppose that the number of nodal domains of a stationary Gaussian field $F:\R^{d}\rightarrow\R$ satisfies the law of large numbers \eqref{eq:NS law large numb}. Then we say that $F$ satisfies the {\em nodal lower concentration property} if, for every $\epsilon>0$,
\begin{equation}
\label{eq:lower conc dec}
\prob\left(\frac{\Nc(F;R)}{\vol B(R)}<c_{NS}(\rho) - \epsilon\right) = o_{R\rightarrow\infty}\left( \frac{1}{R^{d}} \right) .
\end{equation}
\end{definition}

Compared to the law of large numbers~\eqref{eq:NS law large numb}, the nodal lower concentration property \eqref{eq:lower conc dec} {\em quantifies} the decay of the lower tail of $\Nc(F;R)/\vol B(R)$. Rivera--Vanneuville ~\cite[Theorem 1.4]{RiVa} and Beliaev--Muirhead--Rivera ~\cite{BeMuRi} recently proved that $F$ satisfies the nodal lower concentration property provided that the covariance function of $F$ decays sufficiently quickly. In particular, it is sufficient that
\begin{equation*}
r_{F}(x) \le |x|^{-3d-\delta}
\end{equation*}
for some $ \delta>0$ and $x$ sufficiently large.

\begin{theorem}
\label{thm:mean cons vol Eucl}
Let $F:\R^{d}\rightarrow\R$ be a continuous stationary Gaussian field whose spectral measure satisfies $(\rho 1)$--$(\rho 4)$. Let $\Psi=\Psi_{F}$ be the limit volume distribution defined in \eqref{eq:PsiF lim vol mean}. Denote $\Pc=\Pc^{F}$ to be the percolation probability associated to $F$ as in Definition \ref{def:perc prob Euclid}. Then
\begin{enumerate}[(a)]
\item
\label{it:mean vol Eucl lim}
The mean of the limit volume distribution $\Psi_F$ is
\begin{equation}
\label{eq:mean vol lim=(1-perc)/cNS}
\int\limits_{0}^{\infty}(1-\Psi_F(t))dt = (1-\Pc)\cdot \frac{1}{c_{NS}(\rho)}.
\end{equation}
\item
\label{it:empirical mean=1/NS}
If $F$ moreover satisfies the nodal lower concentration property, then the empirical volume mean converges to $\frac{1}{c_{NS}(\rho)}$ in mean, i.e.\
\begin{equation}
\label{eq:empirical mean=1/NS}
\lim\limits_{R\rightarrow\infty}\E\left[\left|\frac{\vol B(R)}{\Nc(F;R)}\cdot\mathbbm{1}_{\Nc(F;R)>0}   -\frac{1}{c_{NS}(\rho)}\right|\right] = 0.
\end{equation}
\end{enumerate}
\end{theorem}

Theorem \ref{thm:mean cons vol Eucl} shows that the first equality of \eqref{eq:mean vol lim Eucl=1/cNS} holds {\em if and only if $\Pc=0$}, much like (the less explicit) Theorem \ref{thm:sum(k mu(k))=2<=>no perc} regarding the connectivity measure. On the other hand, the second equality of \eqref{eq:mean vol lim Eucl=1/cNS} holds for all fields satisfying the nodal lower concentration property regardless of whether $\Pc=0$. We leave open the question of whether in fact the empirical volume mean converges to $1/c_{NS}(\rho)$ in full generality. Indeed it is plausible that all Gaussian fields satisfying $(\rho 1)$--$(\rho 4)$ satisfy the nodal lower concentration property.


\subsection{Ensembles of Gaussian fields on a manifold}

In applications, rather than dealing with a single random field on $\R^d$, one is often given an {\em ensemble} (or sequence) of Gaussian fields, all defined on some fixed Riemannian manifold, that converge to a local limit. In this setting the `escape of mass' for the volume distribution has a slightly different meaning than \eqref{eq:mean vol lim Eucl=1/cNS} (see Theorem \ref{thm:mean cons vol Riem} below).

\smallskip
Let us first make precise the setting in which we work. Let $\M$ be a smooth compact Riemannian $d$-dimensional manifold, and let $\{\Phi_{L}\}_{L\in \Lc}$ be an ensemble of Gaussian fields $\Phi_{L}:\M\rightarrow\R$, with $\Lc\subseteq\R$ some discrete index set. Given a point $x\in\M$ and a sufficiently small neighbourhood $U$ of $x$, we may identify $U$ with a Euclidean sub-domain via the exponential map. More precisely, the exponential map
\[\exp_{x}:T_{x}(\M)\rightarrow \M\]
is a local isometry between a sufficiently small neighbourhood $U\subseteq T_{x}(\M) \equiv \R^{d}$ and $\exp_{x}(U)\subseteq \M$, and by the compactness of $\M$ we can choose $U$ independent of $x$ (under the identification $T_{x}(\M) \equiv \R^{d}$, where we identify $0\in U$ with $x\in\M$). Hence, for every $x$ we may induce a Gaussian field on a domain in $\R^{d}$ and scale it using the linear structure of $\R^{d}$. That is, for $U\subseteq \R^{d}$ so that $\exp_{x}:U\rightarrow \exp_{x}(U)$ is bijective, we define the {\em scaled} Gaussian fields $\Phi_{x;L}:L\cdot U\rightarrow\R$ on the increasing domains
\[ L\cdot U = \{Lu:u\in U\} \]
 to be
\begin{equation}
\label{eq:Phi loc def}
\Phi_{x;L}(u) := \Phi_{L}(\exp_{x}(u/L)).
\end{equation}
The covariance function of $\Phi_{x;L}$ is the function
\begin{equation}
\label{eq:rxL scal covar}
r_{x;L}(u,v):= \E[\Phi_{L}(\exp_{x}(u/L)) \cdot \Phi_{L}(\exp_{x}(v/L))] = r_{\Phi_{L}}(\exp_{x}(u/L),\exp_{x}(v/L)) ,
\end{equation}
defined for $u,v\in L\cdot U$, where $r_{\Phi_{L}}$ is the covariance function of $\Phi_{L}$. Following Nazarov--Sodin \cite{So12,NaSo15}, we consider only the situation in which the ensemble $\{\Phi_{L}\}_{L\in\Lc}$ possesses a `translation invariant local limit':

\begin{definition}[Scaling limits for Gaussian ensembles.]
\label{def:ensemble scal lim}
Let $\{\Phi_{L}\}_{L\in \Lc}$ be a Gaussian ensemble on $\M$, and let $r_{x;L}$ be given by \eqref{eq:rxL scal covar}. We say that $\{\Phi_{L}\}_{L\in \Lc}$ possesses a \textit{translation invariant local limit} as $L\rightarrow\infty$,
if, for almost all $x\in\M$, there exists a continuous covariance kernel $K_{x}:\R^{d}\rightarrow\R$ of a stationary Gaussian field on $\R^{d}$, so that for all $R>0$,
\begin{equation}
\label{eq:loc unif conv covar}
\lim\limits_{L\rightarrow\infty} \sup\limits_{|u|,|v|<R}\left|r_{x;L}(u,v) - K_{x}(u-v)   \right| = 0.
\end{equation}
\end{definition}

Definition \ref{def:ensemble scal lim} is applicable to a number of motivational examples (e.g.\ Kostlan's ensemble, or band-limited fields, see \S\ref{ssec:ens} below). Moreover, in these examples $K_{x}$ is independent of $x$, and so we can associate to the ensemble a single limiting Gaussian field $F:\R^{d}\rightarrow\R$ with covariance $K = K_x$.

\smallskip
Assume now that the above holds (i.e.\ $\{\Phi_{L}\}_{L\in \Lc}$ possesses a translation invariant local limit independent of $x$). Let $\Nc(\Phi_{L})$ denote be the total number of nodal domains of $\Phi_{L}$, and for $t>0$, let $\Nc(\Phi_{L},t)$ denote be the number of those of (Riemannian) volume $<t$. In this setting Nazarov--Sodin ~\cite{So12,NaSo15} proved that
\begin{equation}
\label{eq:NS conv mean man}
\E\left[\left| \frac{\Nc(\Phi_{L})}{L^{d}\vol(\M)}  - c_{NS}(\rho)\right|\right] \rightarrow 0,
\end{equation}
and Beliaev--Wigman ~\cite[Theorem $1.5$]{BeWi} proved that,\footnote{Though stated only for band-limited functions, it is valid in the aforementioned setting.} if $\Psi(\cdot)=\Psi_{F}(\cdot)$ is the cumulative distribution function for $F$ (i.e.\ \eqref{eq:PsiF lim vol mean} is satisfied), then for every continuity point $t$ of $\Psi(\cdot)$, one has
\begin{equation*}
\E\left[\left| \frac{\Nc(\Phi_{L},t/L^{d})}{\Nc(\Phi_{L})}  - \Psi(t) \right|\right] \rightarrow 0,
\end{equation*}
i.e., after the natural scaling, the volume distribution law tends to $\Psi$ in mean.

Since, by the virtue of \eqref{eq:mean vol lim=(1-perc)/cNS} of Theorem \ref{thm:mean cons vol Eucl} applied on $F$, we readily know that
\begin{equation}
\label{eq:1/cNS=perc mean ens}
\frac{1}{c_{NS}(\rho)} = \frac{1}{1-\Pc}\cdot \int\limits_{0}^{\infty}(1-\Psi(t))dt ,
\end{equation}
(i.e.\ the first equality in \eqref{eq:mean vol lim Eucl=1/cNS} holds for the limit law), the question is whether we can relate it to the empirical volume mean
\begin{equation}
\label{eq:emp mean ens}
\frac{\vol(\M)}{\Nc(\Phi_{L})/L^{d}} = \frac{L^{d}\cdot \vol(\M)}{\Nc(\Phi_{L})},
\end{equation}
as asserted in Theorem \ref{thm:mean cons vol Riem} below, for a wide class of ensembles. Note that the corresponding question for mean connectivity trivialises, since we can use Euler's identity on the total nesting graph on $\M$ to verify that the mean connectivity of the nodal domain on $\M$ converges to two in all cases.
Similarly to Theorem \ref{thm:mean cons vol Eucl}, to state our result we need to introduce an analogous `nodal lower concentration property' (c.f.\ Definition \ref{def:nodal low conc}):

\begin{definition}[Nodal lower concentration for ensembles, c.f.\ Definition \ref{def:nodal low conc}]
\label{def:nod low conc ens}
Let $\{\Phi_{L}\}_{L\in \Lc}$ be an ensemble of Gaussian fields possessing a translation invariant local limit as $L\rightarrow\infty$ that is independent of~$x$. Assume further that the spectral measure $\rho$ of the Gaussian field $F$  corresponding to the limit covariance $K$ satisfies $(\rho 1)$--$(\rho 4)$. We say that $\{\Phi_{L}\}_{L\in \Lc}$ satisfies the nodal lower concentration property if, for every $\epsilon>0$,
\begin{equation*}
\prob\left(\frac{\Nc(F;R)}{L^{d}\vol(\M)}<c_{NS}(\rho) - \epsilon\right) = o_{L\rightarrow\infty}\left( \frac{1}{L^{d}} \right).
\end{equation*}
\end{definition}

As an example, it is known~\cite[Theorem $1.1$]{NaSo09} that the random spherical harmonics (see \S\ref{sec:spher harm ARW} below) satisfy the nodal lower concentration property (in fact they satisfy the vastly stronger exponential concentration property). Moreover it was recently shown~\cite{BeMuRi} that the Kostlan ensemble of random homogeneous polynomials (see \S\ref{sec:Kostlan BF ex} below) also satisfies the nodal lower concentration property. On the other hand, imposing the lower concentration property merely on the limit random fields of a Gaussian ensemble (Definition \ref{def:nodal low conc}) is unlikely to yield the lower concentration property for the ensemble, as the former does not control correlations on macroscopic scales.
We are now in a position to state our theorem, asserting the asymptotic equality of \eqref{eq:1/cNS=perc mean ens} and \eqref{eq:emp mean ens} under suitable conditions:

\begin{theorem}
\label{thm:mean cons vol Riem}
Let $\{\Phi_{L}\}_{L\in \Lc}$ be a Gaussian ensemble on $\M$ possessing a translation invariant local limit $K$ as $L\rightarrow\infty$ that is independent of $x$. Assume further that the spectral measure $\rho$ of the Gaussian field $F$  corresponding to the limit covariance $K$ satisfies $(\rho 1)$--$(\rho 4)$, and also that $\{\Phi_{L}\}_{L\in \Lc}$ satisfies the nodal lower concentration property in Definition \ref{def:nod low conc ens}. Then
\begin{equation}
\label{eq:empirical mean=1/NS ens}
\frac{L^d\vol(\M)}{\Nc(\Phi_{L})} \rightarrow \frac{1}{c_{NS}(\rho)} \quad \text{in mean.}
\end{equation}
\end{theorem}

\subsection*{Acknowledgements}

The research leading to these results has received funding from the Engineering \& Physical Sciences Research Council (EPSRC) Fellowship EP/M002896/1 held by Dmitry Beliaev (D.B. \& S.M.), the  EPSRC Grant EP/N009436/1 held by Yan Fyodorov (S.M.), and the European Research Council under the European Union's Seventh Framework Programme (FP7/2007-2013), ERC grant agreement n$^{\text{o}}$ 335141 (I.W.)
We are grateful to P. Sarnak and M. Sodin for the very inspiring and fruitful conversations concerning subjects relevant to this manuscript.


\medskip
\section{Outline of the paper}

In this section we discuss some applications of our main results, and also present an outline of their proofs.

\subsection{Applications}
\label{ssec:ens}
Our results apply to several important ensembles of Gaussian fields on manifolds such as the sphere and torus, as well as to their scaling limits. Some of the applications are rigorous consequences of our main theorems, while others are conjectural.

\subsubsection{Kostlan's ensemble and the Bargmann--Fock limit field}
\label{sec:Kostlan BF ex}
The Kostlan ensemble of degree $n$ homogeneous polynomials is a sequence of Gaussian fields $g_{n}:\R \mathbb{P}^{d}\rightarrow\R$ defined on the real projective space as
\begin{equation}
\label{eq:gn def Eucl}
g_{n}(x) = \sum\limits_{|J|=n}a_{J}\cdot {n\choose J} x^{J},
\end{equation}
where $J=\left(j_{0},j_{1},\ldots,j_{d}\right)$ is a multi-index, $|J|=j_{0}+\ldots+j_{d}$, $x=\left[x_{0}:x_{1}:\ldots :x_{d}\right]$,
$x^{J}=x_{0}^{j_{0}}\cdot \ldots \cdot x_{d}^{j_{d}}$, and $\{a_{J}\}$ are i.i.d.\ standard Gaussians. In the case $d=1$, the study of the zeros of $g_{n}$ is a classical problem in probability theory going back to Shub and Smale ~\cite{ShSm}, and for $d>1$, the study of the nodal structures of $g_{n}$ in the complex algebro-geometric context was initiated by Gayet--Welschinger ~\cite{GaWe}.

\smallskip
Alternatively to \eqref{eq:gn def Eucl}, one may restrict $g_{n}$ on the unit sphere $\Sc^{d}\hookrightarrow \R \mathbb{P}^{d}$, and consider $g_{n}:\Sc^{d}\rightarrow\R$; with this identification, $g_{n}$ is the centred Gaussian field with covariance
\[r_{g_{n}}(x,y):= \E[g(x)\cdot g(y)] = \langle x, y \rangle^{n} = \cos^n(\theta(x,y)),\]
where $\theta(\cdot, \cdot)$ is the angle (or spherical distance) between two spherical points. The upshot is that, with this representation, $g_{n}$ is rotation invariant, with uniformly rapidly decaying correlations, and rapid convergence towards the scaling limit Bargmann--Fock random field $F_{BF}:\R^{d}\rightarrow\R$, with covariance $r_{BF}(x,y)=e^{-\|x-y\|^{2}/2}$. In particular, the ensemble $\{g_{n}\}$ possesses $F_{BF}$ as its translation invariant scaling limit around every point $x\in\Sc^{d}$ (scaling by $\sqrt{n}$). As is evident from its covariance, $F_{BF}$ is stationary and isotropic, with rapid, super-exponential decay of correlations.

\smallskip
In the case $d=2$, it is known ~\cite{BeGa} that the percolation probability $\Pc^{F_{BF}}=0$ of the Bargmann--Fock field vanishes, and, moreover ~\cite{RiVa} the critical level $u^\ast$ is equal to zero. Hence by Theorem \ref{thm:sum(k mu(k))=2<=>no perc} the mean of the limit connectivity measure of $F_{BF}$ and, what is the same, the limit connectivity measure of ${g_{n}}$, are both equal to exactly $2$. In higher dimensions the positivity of $\Pc^{F_{BF}}$ is not known, however, in accordance with Sarnak's insight (explained at the end of \S\ref{ssec:perc}) we believe that $\Pc^{F_{BF}}>0$, so that \eqref{eq:mean lim conn=2} {\em should not} hold. The uniform rapid decay of correlations of both $F_{BF}$ and $\{g_{n}\}$ imply (see the comments after definitions \ref{def:nodal low conc} and \ref{def:nod low conc ens}) the nodal lower concentration property, so that Theorem \ref{thm:mean cons vol Eucl}\eqref{it:empirical mean=1/NS} applies to $F_{BF}$, and Theorem \ref{thm:mean cons vol Riem} applies to $\{g_{n}\}$.

\subsubsection{Spherical harmonics, Arithmetic Random Waves, and their scaling limits}
\label{sec:spher harm ARW}

For $\ell\ge 1$ the degree-$\ell$ spherical harmonics are the harmonic polynomials on $\R^{d+1}$ of degree $\ell$ restricted to the unit sphere $\Sc^{d}$; they constitute a linear space of dimension
\[M_{d;\ell} = \frac{2\ell+d-1}{\ell+d-1} {\ell+d-1 \choose d-1}\]
satisfying the Schr\"{o}dinger equation
\[\Delta_{\Sc^{d}} T_{\ell} + \lambda_{d;\ell} T_{\ell}=0,\]
with (spherical) Laplace eigenvalues $\lambda_{d;\ell}=\ell (\ell + d-2)$. For a $L^{2}$-orthonormal basis $\Ec=\{\eta_{\ell;1},\ldots,\eta_{\ell;M_{d;\ell}}\}$ we define
the random fields on $\Sc^{d}$
\begin{equation*}
T_{\ell}(x) = \frac{1}{\sqrt{M_{d;\ell}}}\sum\limits_{j=1}^{M_{d;\ell}} a_{j}\eta_{\ell;j}(x),
\end{equation*}
with $a_{j}$ standard i.i.d.\ Gaussians; the law of $T_{\ell}$ is independent of the choice of $\Ec$. Equivalently, $T_{\ell}$ is the (uniquely defined) centred Gaussian field on $\Sc^{d}$, with covariance function $\E[T_{\ell}(x)\cdot T_{\ell}(y)] = P_{d;\ell}(\cos(\theta(x,y)))$, where $\theta(x,y)$ is again the spherical distance, and $P_{d;\ell}$ is the degree-$l$ Gegenbauer polynomial (so, in particular, for $d=2$ these are the Legendre polynomials).

\smallskip
The Gaussian ensemble $\{T_{\ell}\}_{\ell\ge 1}$ is important in mathematical physics, cosmology, natural sciences and other disciplines; the fields $T_{\ell}$ appear in the Fourier expansion of any isotropic $L^{2}$-summable Gaussian field on $\Sc^{d}$, hence its importance in the study of the Cosmic Microwave Background (CMB), where $T_{\ell}$, $\ell\rightarrow\infty$, corresponds to high precision experimental measurements. By the standard asymptotics for the Gegenbauer polynomials, $\{ T_{\ell}\}$ possesses a translation invariant local limit, namely, the stationary Gaussian field $F$ with the spectral measure being the hypersurface measure on $\Sc^{d-1}\subseteq\R^{d}$. For example, for $d=2$ these are the planar isotropic monochromatic waves (`Berry's Random Wave Model'), believed ~\cite{Berry77} to represent generic (deterministic) Laplace eigenfunctions on two-dimensional manifolds. For the ensemble $\{T_{\ell}\}$ the exponential nodal concentration was established ~\cite{NaSo09}, stronger than the mere nodal lower concentration property required for the application of Theorem \ref{thm:mean cons vol Riem}. As for the percolation probability, we believe that $\Pc>0$ if and only if $d\ge3$ (see \S\ref{sec:band-lim}).

\smallskip
Another manifold where the solutions for the Schr\"{o}dinger equation can be written explicitly is the $d$-dimensional torus $\Tb^{d}=\R^{d}/[0,1]^{d}$. We may write a general solution to Schr\"{o}dinger equation as
\begin{equation*}
f_{n}(x) = \frac{1}{\sqrt{r_{d}(n)}}\sum\limits_{\|\vec{\lambda}\|^{2}=n} a_{\vec{\lambda}}e(\langle\vec{\lambda},x\rangle),
\end{equation*}
where $n \ge 1$ and the summation is over all lattice points $\vec{\lambda}=(\lambda_{1},\ldots,\lambda_{d})\in\Z^{d}$ satisfying $\|\vec{\lambda}\|^{2}=\lambda_{1}^{2}+\ldots+\lambda_{d}^{2}=n$ (i.e.\ $\vec{\lambda}$ is on a radius $\sqrt{n}$ centred $(d-1)$-hypersphere), $x=(x_{1},\ldots,x_{d})\in\Tb^{d}$, $a_{\vec{\lambda}}\in\C$ are some complex-valued coefficients satisfying
\begin{equation}
\label{eq:a(-lam)=line(a(lam))}
a_{-\vec{\lambda}}=\overline{a_{\lambda}}.
\end{equation}
One can endow $\{f_{n}\}$ with a Gaussian probability measure by taking $\{a_{\vec{\lambda}}\}$ to be standard Gaussian i.i.d.\ (save for \eqref{eq:a(-lam)=line(a(lam))}); the resulting ensemble is referred to as `Arithmetic Random Waves'.

\smallskip
For $d\ge 3$, the ensemble $\{f_{n}(x)\}$ possesses the same translation invariant local limit as $\{T_{\ell}\}$, whereas for $d=2$ this limit arises for generic index sequences, with other scaling limits for exceptional thin index sequences~\cite{Ci,KuWi1}; it is known that the nodal structures of $f_{n}$ are related to the number theoretic properties of these exceptional numbers~\cite{KKW,KuWi2}. For $f_{n}$ the exponential concentration was established by Rozenshein~\cite{Ro} for $d\ge 3$, and $d=2$ with $n$ generic, stronger than needed for an application of
Theorem~\ref{thm:mean cons vol Riem}.

\subsubsection{Band-limited functions}
\label{sec:band-lim}

The examples in \S\ref{sec:spher harm ARW} are particular cases of {\em band-limited} random Gaussian functions for a generic smooth compact $d$-manifold $\M$ (where no spectral degeneracy is expected), put forward by Sarnak--Wigman ~\cite{SaWi}. Let $\Delta$ be the the Laplace--Beltrami operator on $\M$, $\{\varphi_{j}\}_{j\ge 1}$ the (discrete) orthonormal basis of $L^{2}(\M)$ consisting of eigenfunctions satisfying \[\Delta\varphi_{j}+\lambda_{j}\varphi_{j}=0,\]
with corresponding sequence of eigenvalues $\lambda_{j}\ge 0$ nondecreasing, $\lambda_{j}\rightarrow\infty$. Fix a number $\alpha\in [0,1]$ (the `band'), and, given a spectral parameter $T\rightarrow\infty$, we define the $\alpha$-band limited random function to be
\begin{equation}
\label{eq:Phi band-lim def}
\Phi_{T}(x) = \sum\limits_{\alpha\cdot T\le\sqrt{\lambda_{j}}\le T} a_{j}\varphi_{j}(x),
\end{equation}
with $a_{j}$ standard Gaussian i.i.d., where for $\alpha=1$ it is understood that the summation in \eqref{eq:Phi band-lim def} is in the range
\[T- \eta(T) \le \sqrt{\lambda_{j}} \le \sqrt{\lambda_{j}},\]
and $\eta(\lambda)=o(T)$ with $\sqrt{T}\eta(T)\rightarrow\infty$. It is known~\cite{La,Ho,LPS} that $\{\Phi_{T}\}$ possesses a translation invariant local limit, with the limit kernel $K$ being the Fourier transform of the characteristic function of the annulus $$\{y\in\R^{d}:\:  \alpha \le |y|\le 1\},$$ independent of $x$ (for $\alpha=1$, the unit sphere $\Sc^{d-1}\subseteq \R^{d}$). Equivalently, the scaling limit random field $F$ of $\Phi_{T}$ at every point is stationary and rotation invariant (isotropic), and its spectral measure is the characteristic function of the above annulus.

\smallskip
For this fundamental ensemble Sarnak--Wigman ~\cite{SaWi} established a limit connectivity measure $\mu_{\Gamma,d,\alpha}$ on $\Z_{\ge 0}$ that charges all of $\Z_{\ge 0}$, with some extra care required ~\cite{CS} for the case $\alpha=1$ in which the support of the corresponding spectral measure does not contain an interior point; Beliaev--Wigman ~\cite{BeWi} proved the analogous results for the limit volume distribution for nodal domains. 
For the limit random field~$F$, it is not known whether the percolation probability $\Pc$ is positive, nor, in light of the fact that the covariance function $r_{F}$ of $F$ decays too slowly, whether the nodal lower concentration property \eqref{eq:lower conc dec} holds. We believe that the nodal lower concentration property {\em should hold} for both $F$ and $\{\Phi_{T}\}_{T}$ for all $\alpha\in [0,1]$, in all dimensions $d\ge 2$, and, in accordance with Sarnak's insight (explained at the end of \S\ref{ssec:perc}), we believe that $\Pc>0$ if and only if $d\ge 3$, $\alpha\in [0,1]$ arbitrary. If our intuition is correct, the upshot is that \eqref{eq:mean lim conn=2} holds if and only if $d=2$, whereas \eqref{eq:empirical mean=1/NS} and Theorem~\ref{thm:mean cons vol Riem} hold for all $d \ge 2$.

\subsection{Outline of the proofs of the main results}

The proof of Theorem \ref{thm:sum(k mu(k))=2<=>no perc} and Theorem \ref{thm:mean cons vol Eucl}\eqref{it:mean vol Eucl lim} are based on an analysis of the contribution from \textit{boundary components} to, respectively, the total connectivity and volume of the nodal domains. Rather than work with a radius $R$ ball $B(R)$, it will be convenient to redefine $B(R)$ to be the cube $[-R,R]^d$. None of the conclusions in \cite{NaSo15, SaWi, BeWi} are affected by this change; notably each of \eqref{eq:NS conv mean}, \eqref{eq:d(mu(R),mu)>eps->0}, and~\eqref{eq:PsiF lim vol mean} remain valid.

 The proofs of Theorem \ref{thm:sum(k mu(k))=2<=>no perc} and \ref{thm:mean cons vol Eucl} \eqref{it:mean vol Eucl lim} are divided into three steps each:

\smallskip
\noindent \textbf{Step 1.} First we define an appropriate \textit{quantification} of the contribution from the boundary components to the total connectivity and volume; let us focus first on the connectivity. Recall that $G(R) = (V(R), E(R))$ denotes the nesting graph of the nodal domains that are \textit{fully contained} in $B(R)$. We can similarly define the nesting graph $\widebar{G}(R) = (\widebar{V}(R),\widebar{E}(R))$ of \textit{all} nodal domains of the field $F|_{B(R)}$; this is the graph with vertices $\widebar{V}(R)$ the connected components of $B(R)\setminus \Ac(F)$ and edges $\widebar{E}(R)$ which are the connected components of $B(R) \cap \Ac(F)$ that record adjacency among the nodal domains of $F|_{B(R)}$. One advantage of $\widebar{G}(R)$ over $G(R)$ is that it is a.s.\ a tree (by Jordan's Theorem); hence we have by Euler's formula that
\begin{equation}
\label{e:tree}
\sum_{v \in \widebar{V}(R)} \widebar{d}(v) = 2(\widebar{\Nc}(F; R) - 1),
\end{equation}
where $\widebar{\Nc}(F; R) = |\widebar{V}(R)|$ denotes the number of nodal domains of $F|_{B(R)}$.

\begin{definition}[Boundary connectivity]
\label{def:bndr conn}
The boundary connectivity is defined to be
\begin{equation*}
\Ccc(R) :=  \sum_{v \in \widebar{V}(R)} \widebar{d}(v) - \sum_{v \in V(R)}d(v),
\end{equation*}
where $\widebar{d}(v)$ denotes the degree of the vertex $v$ in $\widebar{G}(R)$ (recall that $d(v)$ denotes the degree of the vertex $v$ in $G(R)$).
\end{definition}

Observe that, similarly to \eqref{e:tree}, by Euler's formula
\[  \sum_{v \in V(R)} d(v) = 2(\Nc(F; R) - T(R)) , \]
where $T(R)$ denotes the number of connected components of the union of the closure of all the $D \in V(R)$ (since, in generally, $G(R)$ is a {\em union} of trees). Hence, combining with \eqref{e:tree},
\begin{equation}
\label{eq:C(R)=2(T-1)-}
\Ccc(R) = 2 (T(R)-1) - 2(\widebar{\Nc}(F; R) - \Nc(F; R)) .
\end{equation}
Notice also that $\widebar{\Nc}(F; R) - \Nc(F; R)$ equals the number of nodal domains of $F|_{B(R)}$ that intersect $\partial B(R)$. It is simple to deduce that this has negligible expectation in the limit:

\begin{proposition}
\label{p:boun}
Let $F$ be a continuous stationary Gaussian field with spectral measure $\rho$ satisfying $(\rho 2)$--$(\rho 3)$. Then as $R \to \infty$,
\[  \frac{  \E[ \widebar{\Nc}(F;R) - \Nc(F;R) ] }{\vol B(R)} \to 0  . \]
  \end{proposition}

Together, these observations show that
  \begin{equation}
  \label{e:equivlim}
   \liminf_{R \to \infty}   \frac{\E[\Ccc(R)]}{\vol B(R)}  = \liminf_{R \to \infty} \frac{2\E[T(R)]}{\vol B(R)} \quad \text{and} \quad  \limsup_{R \to \infty}   \frac{\E[\Ccc(R)]}{\vol B(R)}  = \limsup_{R \to \infty} \frac{2\E[T(R)]}{\vol B(R)}
  \end{equation}
 i.e.\ $\E[\Ccc(R)]$ and $2 \E[T(R)]$ are asymptotically equivalent in the large $R$ limit. This fact that will greatly assist the asymptotic analysis of $\Ccc(R)$ that we undertake in Section \ref{s:bc} in order to prove Theorem \ref{thm:mean cons vol Eucl}.

The notion of `boundary volume', analogous to boundary connectivity applied in course of proving Theorem~\ref{thm:sum(k mu(k))=2<=>no perc}, is defined in a significantly simpler manner:

\begin{definition}[Boundary volume]
\label{d:bouvol}
The boundary volume $\Vc(R)$ is the total volume of the connected components of $B(R)\setminus \Ac(F)$ that intersect the boundary $\partial B(R)$.
\end{definition}

Since the nodal set $\Ac(F)$ is a set of zero volume, the boundary volume $\Vc(R)$ can also be expressed as
\[\Vc(R) = \vol( \{ x  \in B(R) : x  \stackrel{F}{\longleftrightarrow} \partial B(R) \} ) . \]
The definitions of boundary connectivity and boundary volume can both be extended to the setting of a Gaussian field $\Phi:\M\rightarrow\R$ on a compact Riemannian manifold $M$, although we formalise this only in the case of the volume:

\begin{definition}[Boundary volume on a manifold]
Fix $x_{0}\in\M$, and $r>0$ sufficiently small so that $\exp_{x_{0}}(\cdot )$ is a bijection on the radius-$r$ ball inside $T_{x_{0}}(\M)$. Then we define $\Vc_{\Phi;x_{0}}(r)$ to be the total volume of the nodal domains of $\Phi$, restricted to the radius-$r$ geodesic ball centred at $x_{0}$, that intersect the boundary of this ball.
\end{definition}

\smallskip
\noindent \textbf{Step 2.} The next step is to link the quantities $\Ccc(R)$ and $\Vc(R)$ to the percolation probability $\Pc$. In the case of the connectivity, the following proposition roughly asserts that the contribution to the connectivity from the boundary is negligible, as a fraction of the total volume of $B(R)$, if and only if $\Pc = 0$:

\begin{proposition}
\label{prop:conbound}
Let $F:\R^{d}\rightarrow \R$ be an a continuous stationary Gaussian field with spectral measure $\rho$ satisfying $(\rho 2)$--$(\rho 3)$ and associated percolation probability $\Pc=\Pc^{F}$. Then
\begin{enumerate}[(a)]
\item
\label{it:conbound1}
If $\Pc=0$, then
\begin{equation*}
\lim\limits_{R\rightarrow\infty} \frac{\E[\Ccc(R)]}{\vol B(R)} = 0.
\end{equation*}
\item
\label{it:conbound2}
Conversely, if $\Pc>0$, and if in addition the spectral measure $\rho$ satisfies $(\rho 4)$, then
\begin{equation*}
\liminf\limits_{R\rightarrow\infty} \frac{\E[\Ccc(R)]}{\vol B(R)} > 0.
\end{equation*}
\end{enumerate}
\end{proposition}

One interpretation of Proposition \ref{prop:conbound} is that if $\Pc > 0$ then the postulated percolating giant nodal domains (see \S\ref{ssec:perc}) make a non-negligible contribution to the total connectivity; we believe it to be of independent interest.

For the volume, we identity a more direct relationship between the contribution from the boundary components and the percolation probability:

\begin{proposition}
\label{prop:volbound}
$\,$
\begin{enumerate}[(a)]
\item
\label{it:volbound1}
Let $F:\R^{d}\rightarrow\R$ be an a continuous stationary Gaussian field with spectral measure $\rho$ satisfying $(\rho 2)$--$(\rho 3)$ and with associated percolation probability $\Pc=\Pc^{F}$. Then
\begin{equation}
\label{eq:bndry vol->perc prob}
\lim\limits_{R\rightarrow\infty} \frac{\E[\Vc(R)]}{\vol B(R)} = \Pc.
\end{equation}
\item
\label{it:volbound2}
Let $\{\Phi_{L}\}_{L\in \Lc}$ be a Gaussian ensemble on $\M$ possessing a translation invariant local limit $K$ as $L\rightarrow\infty$ that is independent of $x$. Suppose the spectral measure $\rho = \rho_F$ corresponding to the limit field $F$ satisfies  $(\rho 2)$--$(\rho 3)$ and has associated percolation probability $\Pc=\Pc^{F}$. Then for every $x_{0}\in \M$,
\begin{equation}
\label{eq:limLlimR perc prob mfld}
\lim\limits_{R\rightarrow\infty} \lim\limits_{L\rightarrow\infty} \frac{\E\left[\Vc_{\Phi_{L};x_{0}}(R/L)\right]}{\vol B(R/L)} = \Pc.
\end{equation}
\end{enumerate}
\end{proposition}

\begin{figure}
\includegraphics[width=0.31\textwidth]{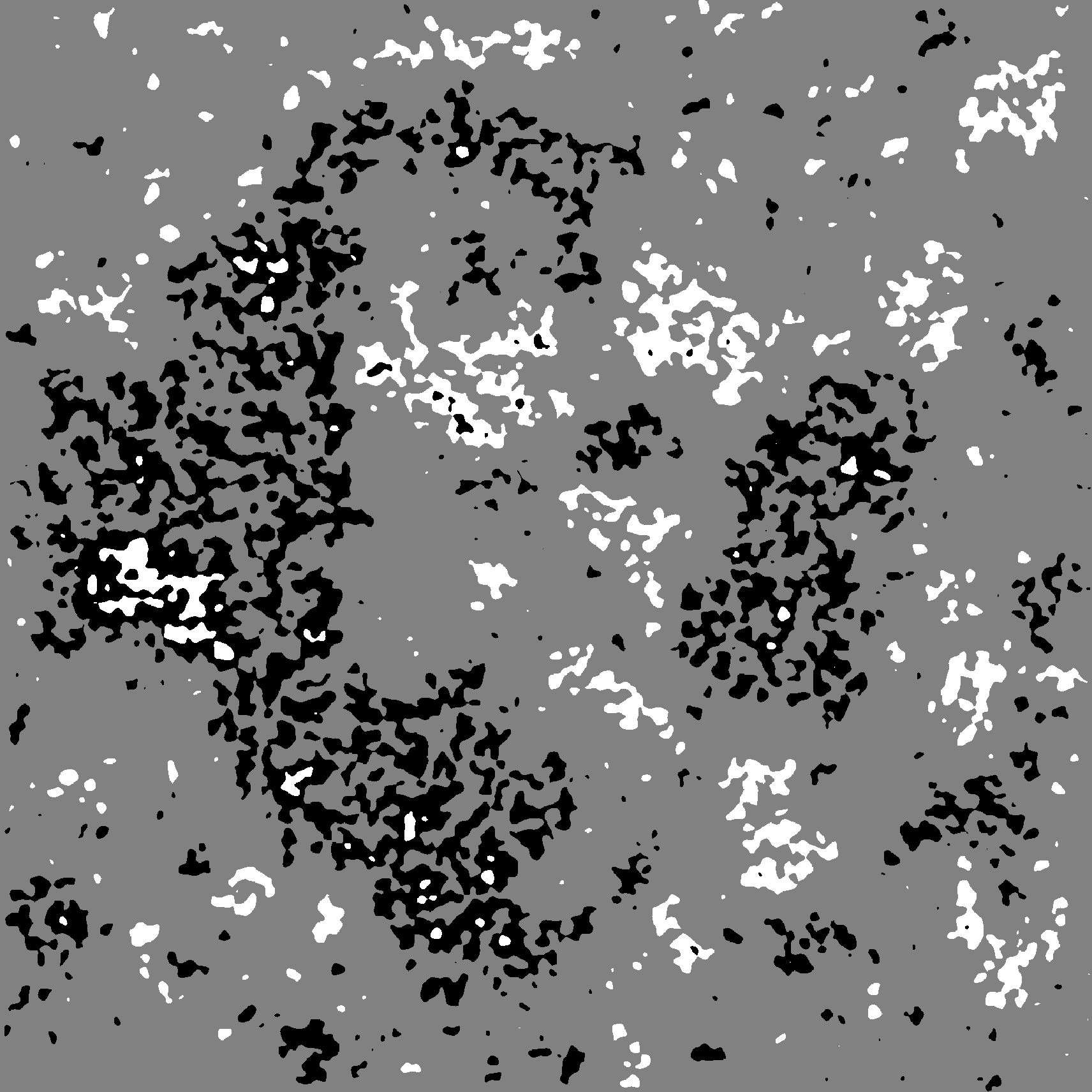}
\includegraphics[width=0.31\textwidth]{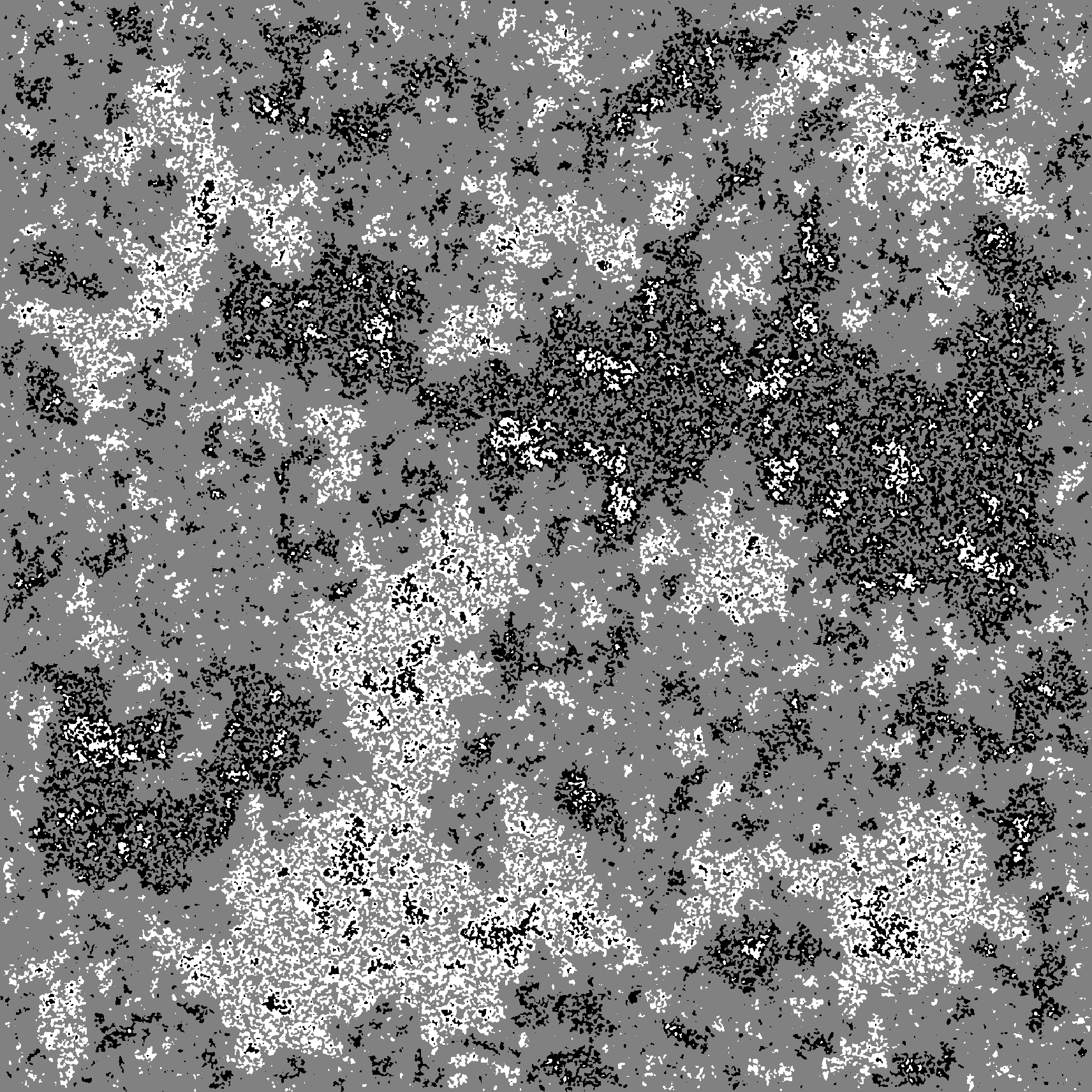}
\includegraphics[width=0.31\textwidth]{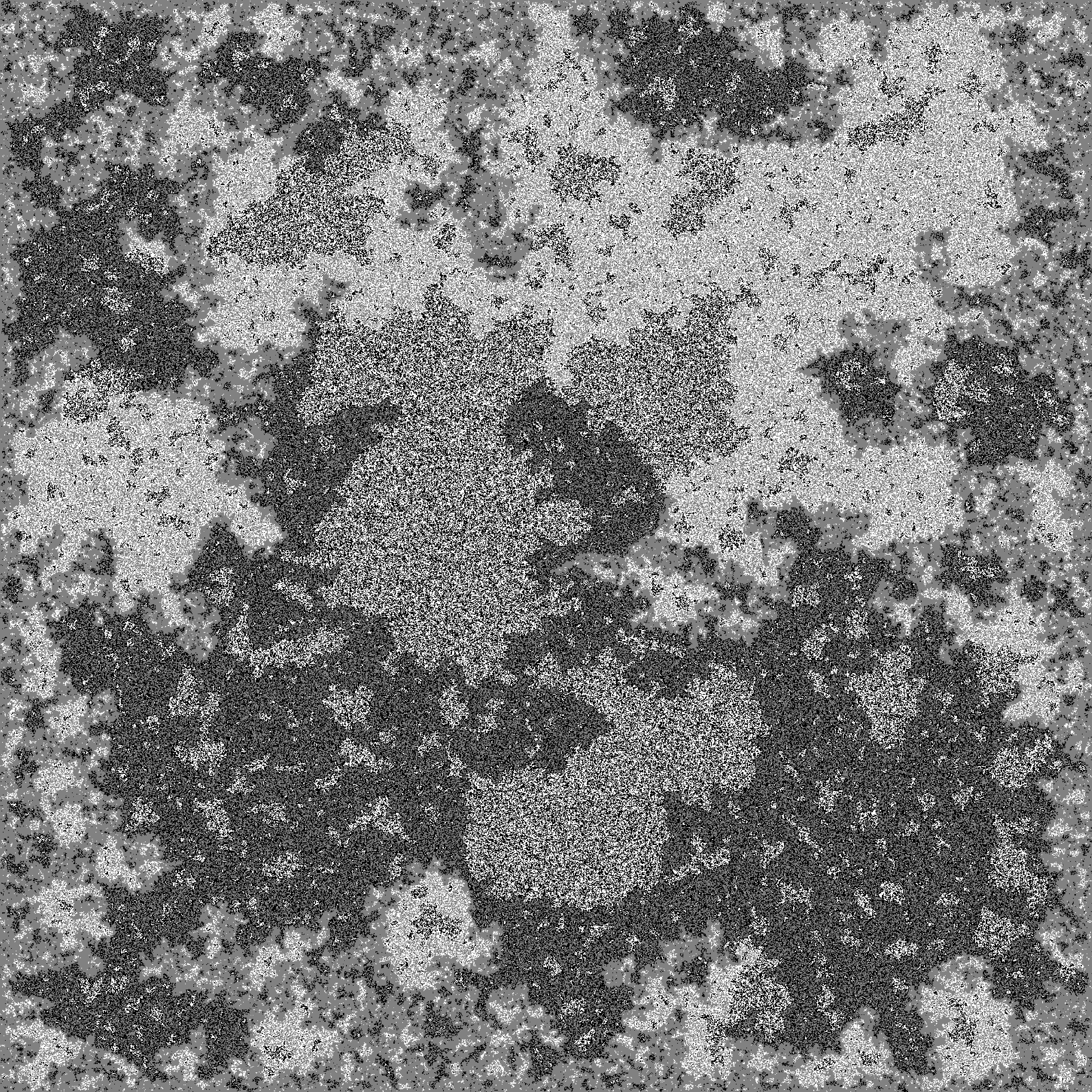}
\caption{Interior nodal domains of the Bargmann-Fock field. The black and white are positive and negative nodal domains respectively, that are entirely contained inside $B(R)$, and the grey are all nodal domains that are connected to the boundary. Note that the grey components are interlaced with most of the black and white domains, including far away from the boundary. Left: $R=200$, middle: $R=800$, right: $R=3200$.  }
\label{fig:boundary}
\end{figure}

Although Proposition \ref{prop:volbound}\eqref{it:volbound2} is not used in the proof of our main theorems, we believe it to be of independent interest in its own right.
Proposition \ref{prop:volbound} implies that in the case $\Pc=0$ the total volume of the nodal components inside $B(R)$ that touch the boundary is negligible. On the other hand, for deterministic reasons there are boundary components of diameter $O(R)$. As illustrated by Figure \ref{fig:boundary}, which shows the boundary components for the Bargmann-Fock field, the typical structure of the boundary components is to have many holes, accounting for their negligible total volume even though their diameter might be large. In particular, the rate of the convergence of the expression on the l.h.s.\ of \eqref{eq:bndry vol->perc prob} (and \eqref{eq:limLlimR perc prob mfld})
to the limit is expected to be slow.

\smallskip
\noindent
\textbf{Step 3.} The final step is to express the mean connectivity and mean volume of the limit measures in terms of the asymptotics formulae for $\Ccc(R)$ and $\Vc(R)$ that appear in Propositions \ref{prop:conbound} and~\ref{prop:volbound} above:

\begin{proposition}
\label{prop:limconvol}
Let $F$ be a continuous stationary Gaussian field with spectral measure $\rho$ satisfying $(\rho 1)$--$(\rho 4)$, and let $\Pc=\Pc^{F}$ be the associated percolation probability. Then
\begin{equation}
\label{eq:lim eq vol}
  \frac{1}{c_{NS}(\rho)} \left(1 -  \limsup\limits_{R\rightarrow\infty}  \frac{\E[\Vc(R)]}{\vol B(R)} \right)  \le \int\limits_{0}^{\infty}\left( 1-\Psi_{F}(t) \right)dt  \le \frac{1}{c_{NS}(\rho)} \left(1 -  \liminf\limits_{R\rightarrow\infty}  \frac{\E[\Vc(R)]}{\vol B(R)} \right)
\end{equation}
and
\begin{equation}
\label{eq:lim eq nest}
 2 -  \frac{1}{c_{NS}(\rho)} \limsup\limits_{R\rightarrow\infty} \frac{\E[\Ccc(R)]}{\vol B(R)}  \le \sum\limits_{k=0}^{\infty}k\cdot \mu_{\Gamma(F)}(k) \le 2 -  \frac{1}{c_{NS}(\rho)} \liminf\limits_{R\rightarrow\infty} \frac{\E[\Ccc(R)]}{\vol B(R)} .
\end{equation}
\end{proposition}

The proof of Propositions \ref{p:boun}, \ref{prop:conbound} and \ref{prop:volbound} will be given in Section \ref{s:bc}, whereas the proof of Proposition~\ref{prop:limconvol} will be given in Section \ref{s:lim}; by combining these propositions we deduce the proofs of Theorem \ref{thm:sum(k mu(k))=2<=>no perc} and Theorem \ref{thm:mean cons vol Eucl}\eqref{it:mean vol Eucl lim}. The proof of Theorem \ref{thm:mean cons vol Eucl}\eqref{it:empirical mean=1/NS} and Theorem \ref{thm:mean cons vol Riem} are more straightforward, and are completed in Section \ref{s:emp}.

\medskip
\section{Analysis of the contribution of boundary components}
\label{s:bc}

In this section we undertake an analysis of the boundary connectivity $\Ccc(R)$ and volume $\Vc(R)$, linking their asymptotics to the percolation probability $\Pc$, and, in particular, prove Propositions \ref{p:boun}, \ref{prop:conbound} and~\ref{prop:volbound}.

\subsection{Proof of Proposition \ref{p:boun}}
The proof of Proposition \ref{p:boun} is standard ~\cite{So12}. The boundary $\partial B(R)$ can be decomposed as the disjoint union of $3^d - 1$ boundary cubes $C_i$ of intermediate dimensions $0 \le i \le d-1$; by standard Morse theory arguments, the number of boundary components
\[  \widebar{\Nc}(F;R) - \Nc(F;R)  \]
is bounded above by the sum, over the boundary cubes $C_i$, of the number of critical points of $F$ restricted to $C_i$. By stationarity and the Kac-Rice formula \cite[Theorem 6.3]{AzWsh} (applicable by $(\rho 2)$--$(\rho 3)$), the expected number of critical points of $F$ restricted to a cube $C$ of dimension $i \ge 1$ is equal to
\[ \text{Vol}(C) \cdot \varphi_{\nabla_C f(0)}(0) \cdot \mathbb{E}[ | \text{det} \nabla_C^2 f(0) \, | \, \nabla_C f(0) =0 ] ,\]
where $\text{Vol}(C)$ is the $i$-dimensional volume of $C$, $\nabla_C$ and $\nabla_C^2$ are respectively the gradient and Hessian of $F$ restricted to $C$, and $\varphi_{\nabla_C f(0)}(0)$ denotes the Gaussian density of $\nabla_C f(0)$ at the value $0$. Since the quantity
\[ s = s(C) = \varphi_{\nabla f(0)}(0) \cdot \mathbb{E}[ | \text{det} \nabla^2 f(0) \, | \, \nabla f(0) =0 ]  \]
depends only on the $i$ directions that span the cube $C$, and since $s(C) > 0$ by $(\rho 2)$--$(\rho 3)$, the expected number of critical points on each boundary cube $C_i$ is proportional to its volume with a constant depending only on the spanning axis directions. Hence in particular
\[  \E[\widebar{\Nc}(F;R) - \Nc(F;R) ] = O(R^{d-1}) .\]

\subsection{Proof of Proposition \ref{prop:volbound}}
The proof of Proposition \ref{prop:volbound} rests on a simple deterministic lemma. For each $R > 0$ and $x \in B(R)$, let $d^-_R(x)$ denote the distance between $x$ and $\partial B(R)$ (i.e. the distance between $x$ and the closest point on $\partial B(R)$ to $x$), and
let $$d^+_R(x) := 2R - d^-_R(x) $$ be the distance between $x$ and the farthest point of $\partial B(R)$, also lying on the axis connecting
$x$ with its closest point on $\partial B(R)$.

\begin{lemma}
\label{l:conv}
Let $g : \R_{\ge 0} \to \R_{\ge 0}$ be a non-increasing function and define $g_\infty = \lim_{s \to \infty} g(s)$. Then, for every $r\in \R$,
as $R\rightarrow\infty$ we have the limits
\begin{equation}
\label{eq:int(g(d-))/B(R)->ginfty}
\frac{  \int_{x \in B(R)} g( (d^-_R(x) - r)^+ )  \, dx  }{\vol B(R)} \to g_\infty
\end{equation}
and
\begin{equation*}
\frac{  \int_{x \in B(R)} g( (d^+_R(x) - r)^+ )  \, dx  }{\vol B(R)} \to g_\infty
\end{equation*}
\end{lemma}

The proof of Lemma \ref{eq:E[Vc(R)]=int perc prob} will be given immediately after the proof of Proposition \ref{prop:volbound}.

\begin{proof}[Proof of Proposition \ref{prop:volbound} assuming Lemma \ref{l:conv}]
Let us begin with part \eqref{it:volbound1}. Observe first that, for each $R > 0$,
\[  \Vc(R) = \vol( \{ x  \in B(R) : x  \stackrel{F}{\longleftrightarrow} \partial B(R) \} )    =  \int_{x \in B(R)} \id \{ x  \stackrel{F}{\longleftrightarrow} \partial B(R) \}   \, dx ,\]
and so
\begin{equation}
\label{eq:E[Vc(R)]=int perc prob}
\E[\Vc(R)]  =   \int_{x \in B(R)} \prob[  x  \stackrel{F}{\longleftrightarrow} \partial B(R) ] \, dx   .
\end{equation}
Notice also that, by the stationarity of $F$, and in light of the fact that for every $x\in B(R)$,
\[  B_{x}(d^{-}_{R}(x))\subseteq B(R) \subseteq B_{x}(d^{+}_{R}(x)),  \]
we have for every $x\in B(R)$,
\begin{equation}
\label{eq:perc prob sand}
\prob[ 0  \stackrel{F}{\longleftrightarrow} \partial B_{d_R^+(x)} ]   \le  \prob[  x \stackrel{F}{\longleftrightarrow} \partial B(R) ] \le  \prob[ 0  \stackrel{F}{\longleftrightarrow} \partial B_{d_R^-(x)} ]  .
\end{equation}
Hence, by substituting \eqref{eq:perc prob sand} into \eqref{eq:E[Vc(R)]=int perc prob}, we obtain the inequality
\begin{equation}
\label{eq:E[V(R)]/B(R) sandwich}
\frac{  \int_{x \in B(R)}  \prob[ 0  \stackrel{F}{\longleftrightarrow}  \partial B_{d_R^+(x)}]    \, dx  }{ \vol B(R)  }   \le   \frac{\E[\Vc(R)]}{\vol B(R)}  \le       \frac{  \int_{x \in B(R)}  \prob[ 0  \stackrel{F}{\longleftrightarrow}  \partial B_{d_R^-(x)} ]   \, dx  }{ \vol B(R)  }  .
\end{equation}
Applying Lemma \ref{l:conv} to the non-increasing function $g(s) =   \prob[  0 \stackrel{F}{\longleftrightarrow} \partial B(s) ]$, yields that both the l.h.s.\ and the r.h.s.\ of \eqref{eq:E[V(R)]/B(R) sandwich} converge, as $R \to \infty$, to the limit
\[ \Pc := \lim_{s \to \infty}  \prob[  0 \stackrel{F}{\longleftrightarrow} \partial B(s) ], \]
and therefore so does $\E[\Vc(R)]/\vol B(R)$, which is the statement of Proposition \ref{prop:volbound}~\eqref{it:volbound1}.

\smallskip
We now turn to part \eqref{it:volbound2}. Recall the definition of the scaled random fields $\Phi_{x;L}$ in \eqref{eq:Phi loc def}, with covariance $r_{x; L}$ given by \eqref{eq:rxL scal covar}.
The assumed locally-uniform convergence \eqref{eq:loc unif conv covar} of the covariance kernels $r_{x; L}$ to $K$ ensure that, on any compact domain, the random field $\Phi_{x;L}$ converges in law, in the $C^0$ topology, to the translation invariant local limit field $F$. Next we observe that the function $h$ that maps a $C^3$-smooth function on a piece-wise smooth compact domain $D \subseteq \mathbb{R}^d$ to the total volume of the nodal domains that intersect $\partial D$ is continuous in the $C^0$ topology up to a null set of~$F$. This is since the set of discontinuities of $h$ is contained in the set of functions such that there is a critical point of $F|_D$ or $F|_{\partial D}$ with height zero, which is indeed a null set for~$F$ by Bulinskaya's Lemma, valid by $(\rho 2)$--$(\rho 3)$ (see e.g.~\cite[Proposition 6.12]{AzWsh}).

Hence, by the Continuous Mapping Theorem, we have the convergence in law
\begin{equation}
\label{eq:vol fin->law vol lim}
\frac{\Vc_{\Phi_L;x_0}(R/L)}{\vol B(R/L)}\stackrel{\text{Law}}{\longrightarrow} \frac{\Vc(R)}{\vol B(R)}.
\end{equation}
Since the random variables on both r.h.s.\ and l.h.s.\ of \eqref{eq:vol fin->law vol lim} are clearly bounded, their means also converge, i.e.\
\[ \lim_{L \to \infty} \frac{\E[\Vc_{\Phi_L;x_0}(R/L)]}{\vol B(R/L)}  = \frac{\E[\Vc(R)]}{\vol B(R)} . \]
In light of part \eqref{it:volbound1} of Proposition \ref{prop:volbound}, we have the result upon taking the limit $R \to \infty$.
\end{proof}

\begin{proof}[Proof of Lemma \ref{l:conv}]
Since $g$ is non-increasing, we have the inequality
\begin{equation*}
\frac{  \int_{x \in B(R)} g( (d^\pm_R(x) - r)^+ )  \, dx  }{\vol B(R)} \ge g_{\infty}.
\end{equation*}
Hence in light of the trivial inequality $d^{+}(\cdot) \ge d^{-}(\cdot)$, to prove both statements of Lemma \ref{l:conv}, it is sufficient to prove \eqref{eq:int(g(d-))/B(R)->ginfty} for $r>0$ only. Moreover, without loss of generality we may assume that $g(0)=1$ and $g_\infty = 0$ (as otherwise we may pass to $(g(\cdot)-g_{\infty})/(g(0)-g_\infty)$). Integrating over $d-1$-dimensional cubic shells (technically justified by dividing the cube into $2d$ identical right-pyramids and applying the smooth co-area formula),
\begin{align*}
\int_{x \in B(R)} g( (d^-_R(x) - r)^+ ) \, dx & = 2^d d \int_{s = 0}^{R-r} s^{d-1} g((R - r - s )  \, ds  + 2^d d \int_{s = R-r}^{R} s^{d-1}   \, ds  \\
 & = 2^d d \int_{s = 0}^{R-r} s^{d-1} g( R - r - s)  \, ds  + 2^d ( R^d - (R-r)^d) ,
\end{align*}
and so it remains to show that for every $\varepsilon, r > 0$ there exists an $R > 0$ sufficiently large so that
\[ \int_{s = 0}^{R-r}  s^{d-1} g(R - r - s) \,  ds  <   \varepsilon  R^d .\]
Let $t > 0$ be such $g(t) \le (d/2) \varepsilon$, and recall that $g$ is bounded by $1$. Then
\begin{align*}
 \int_{s = 0}^{R-r} s^{d-1} g(R - r - s)  ds & =   \int_{s = 0}^{R-r-t} s^{d-1} g(R - r - s)  \, ds + \int_{s = R - r - t}^{R-r} s^{d-1} g(R - r - s) \, ds \\
 & \le (d/2) \varepsilon \int_{s = 0}^{R-r-t} s^{d-1} \, ds +  \int_{s = R - r - t}^{R-r}  s^{d-1}  \, ds \\
  & \le (d/2) \varepsilon \int_{s = 0}^{R} s^{d-1} \, ds +  \int_{s = R - r - t}^{R-r}  R^{d-1}  \, ds \\
 & \le (\varepsilon/2)  R^d + t R^{d-1},
 \end{align*}
which is less than $ \varepsilon  R^d$ for sufficiently large $R$.
\end{proof}

\subsection{Proof of Proposition \ref{prop:conbound}}

As a preparation towards proving Proposition \ref{prop:conbound} we will need three auxiliary lemmas. The first is a simple deterministic bound on the connected components of a set $S \subseteq B(R)$:

\begin{lemma}
\label{l:squares}
There exists an absolute constant $c > 0$ with the following property.
For every $R > 1$ and $S \subseteq B(R)$ closed subset, and $\varepsilon > 0$, the number of connected components of $B(R) \setminus S$ whose volume is at least $\varepsilon$ is bounded above by
\[   K(S)\cdot ( \varepsilon^{-1} + 1) + c R^{d-1} ,   \]
where
\begin{equation}
\label{eq:K(S) intersect def}
K(S) := \# \left\{ k \in \mathbb{Z}^d : S  \text{ intersects the cube } [0,1]^d + k   \right\}  .
\end{equation}
\end{lemma}

Our next lemma, borrowed from \cite{So12,SaWi}, shows that, under the usual assumptions on~$F$, the limit volume distribution $\Psi_F$ exhibits at most power behaviour at the neighbourhood of the origin, i.e.\ yields an upper bound for the (asymptotic) number of {\em small nodal domains}:

\begin{lemma}
\label{l:small}
Let $F$ be a continuous stationary Gaussian field with spectral measure $\rho$ satisfying $(\rho 2)$--$(\rho 3)$. Then there exist constants $c_1, c_2, t_0 > 0$ such that, for all $t_0 > t > 0$,
\[  \Psi_F(t) < c_1 t^{c_2}  .\]
\end{lemma}

Finally we state a simple consequence of $\Pc > 0$, namely that it guarantees the existence, with positive probability, of a nodal domain that (i) lies fully inside a small ball $B(r)$ and (ii) is connected to the boundary $\partial B(R)$ by another nodal domain:

\begin{lemma}
\label{l:w}
Assume that $\Pc > 0$, and assume also that the spectral measure $\rho$ satisfies $(\rho 2)$--$(\rho 4)$. Recall that $V(r)$ denotes the set of nodal domains that are fully contained within $B(r)$. Then there exists a number $r > 0$ such that
\[ \liminf_{R \to \infty} \, \prob \big[ \exists D \in V(r) :  \partial D \stackrel{F}{\longleftrightarrow} \partial B(R) \big] > 0 . \]
\end{lemma}

We are now ready to prove Proposition \ref{prop:conbound}. Recall that $V(R)$ is the set of nodal domains that are fully contained within $B(R)$, and $\widebar{V}(R)$ is the set of nodal domains of the field $F$ restricted to $B(R)$; hence $\widebar{V}(R)\setminus V(R)$ is the set of nodal domain of $F|_{B(R)}$ that intersect $\partial B(R)$. Recall also that $\Nc(F;R) = |V(R)|$ and $\widebar{\Nc}(F;R) = |\widebar{V}(R)|$.

\begin{proof}[Proof of Proposition \ref{prop:conbound} assuming lemmas \ref{l:squares}-\ref{l:w}]

We begin with part \eqref{it:conbound1}. For each $R > 0$, define
\[   S_R := \cup_{D_i \in \widebar{V}(R) \setminus V(R)}  D_i   = \{ x \in B(R) : x \stackrel{F}{\longleftrightarrow} \partial B(R) \} \]
to be the union of the nodal domains  of $F$ restricted to $B(R)$ that intersect the boundary. Let $T(R)$ be the number of connected components of $B(R) \setminus S_R$, and observe that this agrees with the definition given immediately after Definition \ref{def:bndr conn}. Equation \eqref{eq:C(R)=2(T-1)-} states that
\begin{equation*}
 \Ccc(R) = 2 (T(R)-1) - 2(\widebar{\Nc}(F; R) - \Nc(F; R))
 \end{equation*}
and hence, in view of Proposition \ref{p:boun}, to establish the result it is sufficient to show that
\begin{equation}
\label{e:tr}
 \lim_{R \to \infty} \frac{ \E [ T(R) ] }{ \vol B(R) } = 0 .
 \end{equation}
Applying Lemma \ref{l:squares}, for every $\varepsilon > 0$ we have that
\begin{equation}
\label{eq:T(R) bnd N, K(SR)}
T(R) \le   \Nc(F, \varepsilon; R)   +   K(S_{R}) ( \varepsilon^{-1} + 1) + c R^{d-1} ,
\end{equation}
where $K(S_{R})$ was introduced in \eqref{eq:K(S) intersect def},
and $c > 0$ is an absolute constant. Suppose that $\varepsilon > 0$ is a continuity point of the limit volume distribution~$\Psi_F$. By \eqref{eq:PsiF lim vol mean}, as $R \to \infty$,
\[   \frac{ \E[  \Nc(F, \varepsilon; R]}{  \vol B(R) } \to  c_{NS}(\rho)\cdot \Psi_F(\varepsilon).  \]
Since  Lemma \ref{l:small} implies that, as $\varepsilon \to 0$,
\[ \Psi_F(\varepsilon) \to 0,    \]
we deduce that
\begin{equation}
\label{e:nbou}
 \lim_{\varepsilon \to 0 } \lim_{R \to \infty} \frac{ \E[ \Nc(F, \varepsilon; R) ] +  c R^{d-1}}{ \vol B(R)  }   = 0 ,
 \end{equation}
where the limit as $\varepsilon \to 0$ is understood as being taken on a subsequence of continuity points of~$\Psi_F$.

Turning to bounding $K(S_{R})$, we first claim that if $\Pc = 0$ then
\[   \lim_{R \to \infty}  \prob[ D \stackrel{F}{\longleftrightarrow} \partial B(R) ]  = 0  \]
for an arbitrary compact domain $D$. To this end, we observe that since $0$ does not lie on a nodal component a.s., the nodal domain containing $0$ covers a small cube $B(\varepsilon)$ with probability tending to $1$ as $\varepsilon \to 0$. Hence if $\Pc = 0$ then it cannot be the case that
\[   \liminf_{R \to \infty}  \prob[ B(\varepsilon)  \stackrel{F}{\longleftrightarrow} \partial B(R) ]  > 0 , \]
for arbitrary small $\varepsilon > 0$, since then
\[   \liminf_{R \to \infty}  \prob [ 0  \stackrel{F}{\longleftrightarrow} \partial B(R) ]  > 0 , \]
 which is in contradiction with $\Pc = 0$. Thus we have that
\[   \lim_{R \to \infty}  \prob[ B(\varepsilon)  \stackrel{F}{\longleftrightarrow} \partial B(R) ]  = 0  \]
for sufficiently small $\varepsilon > 0$, and we deduce the claim by covering $D$ with a finite number $C_i$ of translations of $B(\varepsilon)$.

Next, applying Lemma \ref{l:conv} to the function $g(s) =  \prob[   B(1)  \stackrel{F}{\longleftrightarrow}  \partial B(s) ]$ and the constant $r = 1$, and arguing as in the proof of Proposition \ref{prop:volbound}, we deduce that
\begin{equation}
\label{e:c}
 \lim_{R \to \infty} \frac{\E[K(S_{R})]}{\vol B(R)} = 0 .
 \end{equation}
Combining \eqref{e:nbou} and \eqref{e:c} and substituting these into \eqref{eq:T(R) bnd N, K(SR)},
while sending first $R \to \infty$ and then $\varepsilon \to 0$, we arrive at \eqref{e:tr}.

\vspace{2mm}

Let us now establish part \eqref{it:conbound2}. As in part \eqref{it:conbound1}, it is sufficient to prove that
\begin{equation}
\label{e:cb}
\liminf_{R \to \infty} \frac{\E[ T(R) ]  }{\vol B(R)}  > 0 .
\end{equation}
Fix $r > 0$ as in Lemma \ref{l:w} and consider tiling $B(R)$ with $O(R^d)$ disjoint translations $E_i$ of $B(r)$ (ignoring the leftover untiled space).
Observe that, since $E_i$ are disjoint,
\[ \E[T(R)]  \ge \sum_{E_i}  \prob[ \exists D \subseteq V(E_i) :   \partial D  \stackrel{F}{\longleftrightarrow}  \partial B(R) ] . \]
Hence applying Lemma \ref{l:conv} to the function
\[ g(s) =   \prob[   \exists D \in V(E_i) :   \partial D  \stackrel{F}{\longleftrightarrow}  \partial B(s) ]  \]
 and the constant $r$, and arguing as in the proof of Proposition \ref{prop:volbound}, we deduce that
\begin{equation}
\label{eq:limE(C(R))/vol(B(R))>=ginf}
 \liminf_{R \to \infty} \frac{\E[ T(R) ]  }{\vol B(R)} \ge  \lim_{s \to \infty} g(s)  .
\end{equation}
Since $\lim_{s \to \infty} g(s)$ is strictly positive by the virtue of Lemma \ref{l:w}, so is the l.h.s.\ of \eqref{eq:limE(C(R))/vol(B(R))>=ginf}, yielding \eqref{e:cb}.
\end{proof}

We now prove the auxiliary Lemmas \ref{l:squares}-\ref{l:w}:

\begin{proof}[Proof of Lemma \ref{l:squares}]
Let $\Kc(S)$ denote the union of all cubes $ [0,1]^d + k $, $k\in\Z^{d}$, that intersect $S$;
by the definition \eqref{eq:K(S) intersect def} of $K(\cdot)$
we have
\begin{equation}
\label{eq:K(S)=vol(Kc(S))}
K(S)=\vol(\Kc(S)).
\end{equation}
A connected component of $B(R) \setminus S$ is either contained within $\Kc(S)$ or is not. The number of components contained within
$\Kc(S)$ with volume at least $\varepsilon$ is bound above by
\begin{equation}
\label{eq:bnd contained comp}
\vol(\Kc(S)) \cdot \varepsilon^{-1} =  K(S)  \cdot \varepsilon^{-1} ,
\end{equation}
by \eqref{eq:K(S)=vol(Kc(S))}.

On the other hand, we may bound the number of those components not lying inside $\Kc(S)$ by invoking the geometric Lemma \ref{lem:comp monotone} below (with $S$ taking the role of $A$ and $\Kc(S)$ taking the role of $B$) to be at most $K(S) + c R^{d-1}$. Together with the bound \eqref{eq:bnd contained comp} for those components lying inside $\Kc(S)$,
this yields the statement of Lemma \ref{l:squares}.
\end{proof}

\begin{lemma}
\label{lem:comp monotone}
Let $B\subseteq B(R)$ be a finite union
\[
B=\bigcup\limits_{k\in I}(k+[0,1]^{d})
\]
of cubes of the form $k+[0,1]^{d}$,
$I\subseteq \Z^{d}$, and $A\subseteq B$ a closed set.
Then the number of connected components of $B(R) \setminus A$ not contained within $B$
is at most
\begin{equation}
\label{eq:vol(B)+cR^d-1}
\vol(B) + c\cdot R^{d-1}.
\end{equation}
\end{lemma}

\begin{figure}
\includegraphics[width=0.4\textwidth]{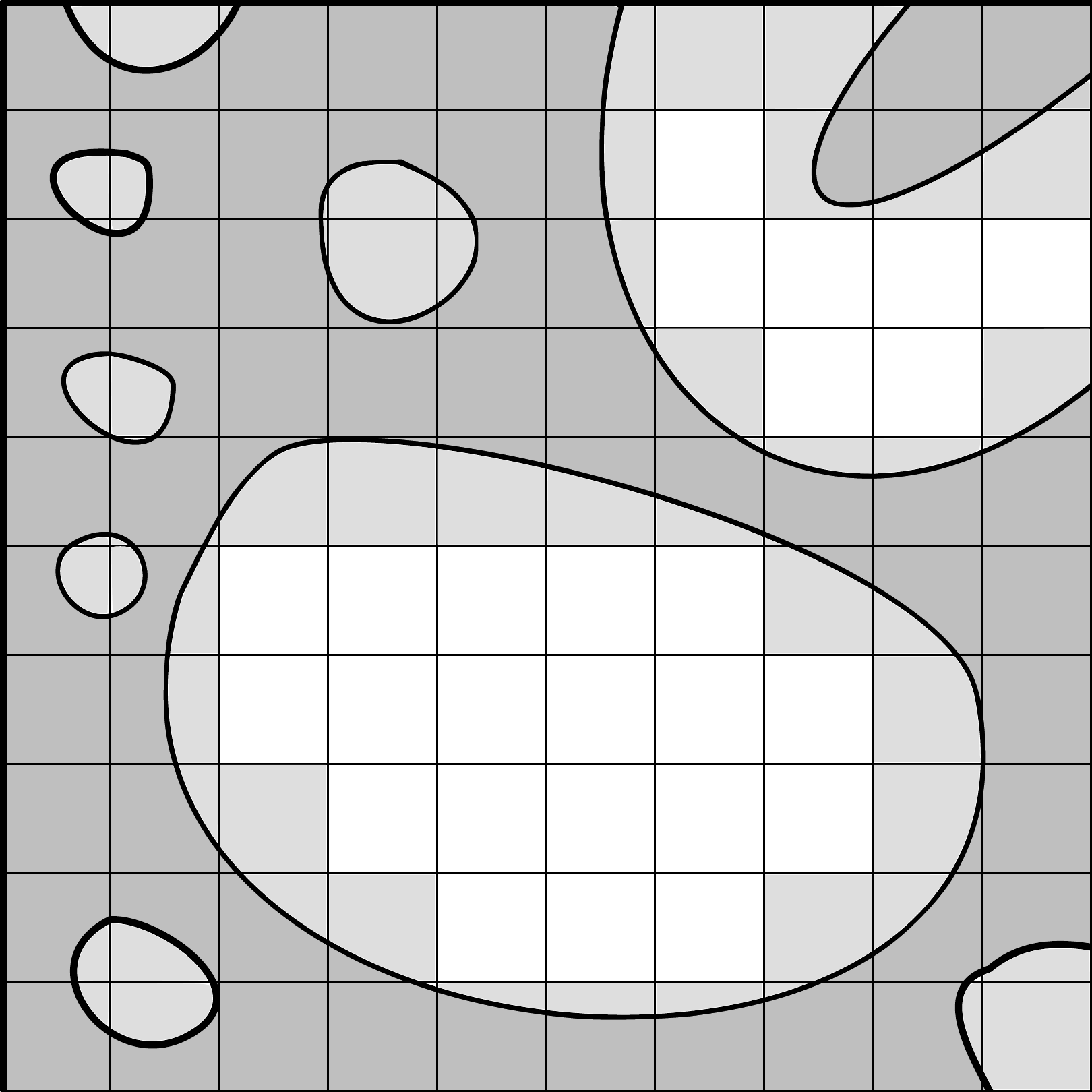}
\caption{An illustration of Lemma \ref{lem:comp monotone}, depicting the set $B$ (in grey) and the set $A \subseteq B$ (in dark grey). The number of connected components of $B(R)\setminus A$ that are not contained within $B$ can be bounded by the number of connected components of $B(R) \setminus B$ (in white), which in turn can be bounded (up to boundary effects) by the volume of $B$; in this figure there are two such components.}
\label{fig:lemma 3.2}
\end{figure}

\begin{proof}
First, every connected component of $\widetilde{A}:=B(R)\setminus A$ that is not contained within $B$ intersects $\widetilde{B}:=B(R)\setminus B$ (as otherwise it would be fully contained in $B$), and hence contains a distinct connected component of $B(R)\setminus B$. This induces an injection between the connected components of $\widetilde{A}$ and the connected components of $\widetilde{B}$, so the number of the former is bounded from above by the latter. Now we claim that the bound \eqref{eq:vol(B)+cR^d-1} is applicable for the number of connected components of $\widetilde{B}$.

To this end we associate to each connected component $C$ of $\widetilde{B}$ a distinct cube that is either in~$B$ or adjacent to $\partial B(R)$ in the following manner. Given a point $x=(x_{1},x_{2},\ldots, x_{d})\in C$ we increase $x_{1}$ until we escape $C$, i.e.\ find the smallest $x_{1}'>x_{1}$ so that $x'=(x_{1}',x_{2},\ldots,x_{d})\notin C$; since $C$ is open (being a complement of a closed set), $x'$ is in the
interior of one of the faces of the cubes $[0,1]^{d}+k\subseteq B$, or one of the cubes $[0,1]^d+k$ intersecting $\partial B(R)$. These cubes are clearly distinct for different components $A$, so their number is bounded by \eqref{eq:vol(B)+cR^d-1}, as claimed.
\end{proof}

\begin{proof}[Proof of Lemma \ref{l:small}]
This is a restatement of \cite[Lemma 9]{So12} (the full proof given in \cite[Lemma 4.12]{SaWi}). Although the result in \cite{SaWi} is stated only for certain special cases of $F$, its proof holds unimpaired for all $F$ satisfying the axioms $(\rho2)$--$(\rho3)$.
It also yields the universality of the exponent $c_2$, depending only on the dimension $d$ (and the threshold $t_0$), although the constant $c_1$ also depends on the field $F$.
\end{proof}

\begin{proof}[Proof of Lemma \ref{l:w}]
\begin{figure}
\includegraphics[width=0.4\textwidth]{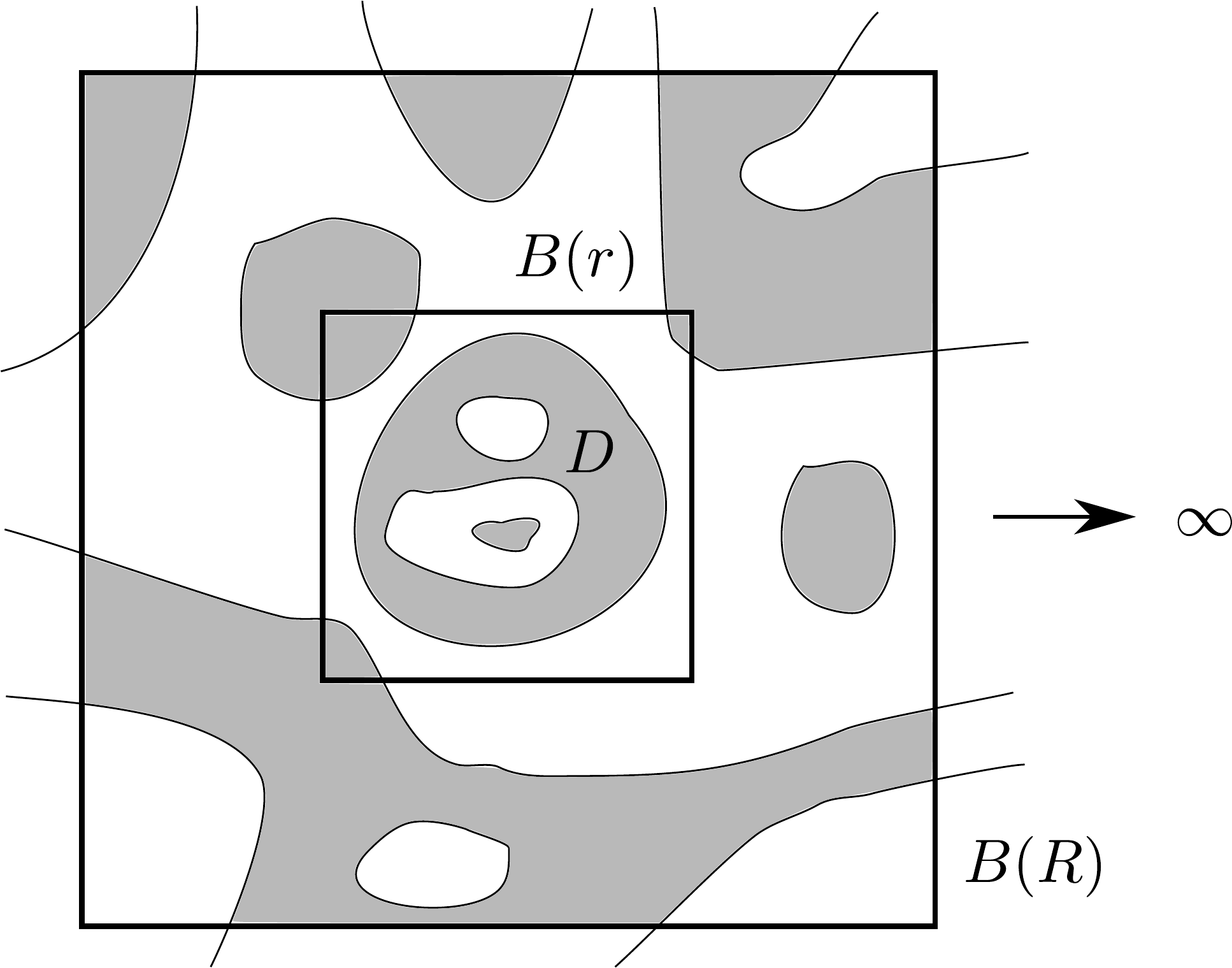}
\caption{With positive probability there is a nodal domain $D$ lying entirely inside $B(r)$, and which is connected by an infinite component to the boundary of an arbitrary large square $B(R)$. }
\label{fig:lemma 3.4}
\end{figure}

Let $\Ec$ denote the event that the positive excursion set
\[ \{ x \in \R^d : F(x) > 0 \} \]
 has an unbounded connected component. Since the percolation probability $\Pc$ is positive, and by the symmetry of $F$ and $-F$, the event $\Ec$ has positive probability. Since moreovoer the random field $F$ is ergodic (equivalent to $(\rho 1)$), and the event $\mathcal{E}$ is translation invariant, we deduce that $\Ec$ occurs a.s.\

Now, let $S$ denote the union of all the unbounded connected components of the positive excursion set $\{ x \in \R^d : F(x) > 0 \}$. Since we assumed $(\rho 4)$ that $c_{NS}(\rho)>0$, there exists a positive density of bounded nodal domains, and so $0 \notin S$ with positive probability. Hence, for sufficiently large $r > 0$, and letting $W$ denote the component of $S^c$ that contains $0$, the event
\[
\Fc = \{0 \in S^c \text{ and } W \subset B(r)\}
\]
holds with positive probability. The set $W$ is the union of the nodal domain $D$ with all nodal domains that are inside of $D$, see Figure \ref{fig:lemma 3.4}.

Finally, assume the event $\Ec \cap \Fc$, and notice that $W$ contains a nodal domain in $V(r)$ that has $\partial W$ as a boundary component. Since $W$ is a component of $S^c$, it must be the case that $\partial W \stackrel{F}{\longleftrightarrow} \infty$, and so
\[ \Ec \cap \Fc \subseteq \cap_{R > r} \{ \exists D \in V(r) :  \partial D \stackrel{F}{\longleftrightarrow} \partial B(R)  \}  . \]
Since $\prob(\Ec \cap \Fc) = \prob(\Fc) > 0$, we deduce the result.
\end{proof}

\bigskip
\section{The mean connectivity and volume of the limit distribution}
\label{s:lim}

In this section we show how to express the mean of the limit connectivity and volume measures $\mu_{\Gamma(F)}$ and $\Psi_F$ in terms of the asymptotic formulae for $\Ccc(R)$ and $\Vc(R)$ (that appear in propositions \ref{prop:conbound} and~\ref{prop:volbound}); in particular, we prove Proposition \ref{prop:limconvol}.
Recall that, by \eqref{eq:d(mu(R),mu)>eps->0}, for every $k\ge  0$,
\begin{equation*}
\mu_{\Gamma(F);R}(k) = \frac{1}{|V(R)|}\cdot \# \{ v\in V(R) : d(v) = k \} \to  \mu_{\Gamma(F)}(k)
\end{equation*}
in probability, as $R \to \infty$. Since we know from \eqref{eq:NS conv mean} that
\begin{equation}
\label{e:llnvol2}
 \frac{|V(R)|}{\vol B(R)} = \frac{\Nc(F;R)}{\vol B(R)}  \to c_{NS}(\rho) \quad \text{in mean,}
 \end{equation}
by the triangle inequality we can deduce, similarly to \eqref{eq:PsiF lim vol mean}, that
\begin{equation}
\label{e:llncon2}
\frac{ \E[  \# \{ v\in V(R) : d(v) = k \}] }{ c_{NS}(\rho) \cdot \vol B(R) }  \to  \mu_{\Gamma(F)}(k)
\end{equation}
We restate \eqref{eq:PsiF lim vol mean} for convenience in the form
\begin{equation*}
\frac{ \E[  \# \{ v\in V(R) : \vol(v) <t  \}] }{ c_{NS}(\rho) \cdot \vol B(R) }  \to  \Psi_F(t),
\end{equation*}
valid at all continuity points of $\Psi_F$; equivalently, in light of \eqref{e:llnvol2}
\begin{align}
\label{e:lim vol>t}
\frac{ \E[  \# \{ v\in V(R) : \vol(v) \ge t  \}] }{ c_{NS}(\rho) \cdot \vol B(R) }  & = \frac{ \E[ \Nc(F; R)] - \E[  \# \{ v\in V(R) : \vol(v) < t  \}] }{ c_{NS}(\rho) \cdot \vol B(R) } \\
\nonumber &\to  1-\Psi_F(t).
\end{align}
We remark that passing to the complement \eqref{e:lim vol>t} (i.e.\ working with domains of volume $\ge t$ and not $<t$) is an important technical step since it will eventually allow us to invoke the Monotone Convergence Theorem when working with the convergent integral $\int_0^\infty(1-\Psi_{F}(t)) dt$ in the proof of Proposition \ref{prop:limconvol} below (see \eqref{e:limits1}).

\begin{proof}[Proof of Proposition \ref{prop:limconvol}]
We first prove statement \eqref{eq:lim eq vol} of Proposition \ref{prop:limconvol}.
To this end we let
\[R_i = 2^i,\]
and partition the cube
\[B(R_i):=[-R_{i},R_{i}]^{d}\]
 into $2^d$ disjoint cubes $C_{i-1;j}$, $j=0,\ldots, 2^{d}-1$ of side length $2R_{i-1}$. We extend the notation of the nesting graph $G(R) = (V(R), E(R))$ to cover the cubes $C_{i-1;j}$, i.e.\ define $$G(C_{i-1;j}) = (V(C_{i-1;j}), E(C_{i-1;j}))$$ analogously to $G(R) = (V(R), E(R))$. Notice that every nodal domain that is fully inside $B(R_i)$ is either fully inside one of the $C_{i-1;j}$ or intersects the boundary of at least one of the $C_{i-1;j}$. Hence, neglecting the latter domains, we have for every $t>0$,
\begin{equation}
\label{eq:vol v sup conv}
\# \{ v\in V(R_i) : \vol(v) \ge t \} \ge \sum\limits_{j=0}^{2^{d}-1} \# \{ v\in V(C_{i-1;j}) : \vol(v) \ge t \} .
\end{equation}

Taking expectations of both sides of \eqref{eq:vol v sup conv}, and upon exploiting the stationarity of $F$, this implies that
\[ \E[\# \{ v\in V(R_i) : \vol(v) \ge t \} ] \ge 2^d \E[\# \{ v\in V(R_{i-1}) : \vol(v) \ge t \} ], \]
which in turn implies that the sequence
\begin{equation*}
\varphi^i_F(t):= \frac{\E[\# \{ v\in V(R_i) : \vol(v) \ge t \}] }{ c_{NS}(\rho)\cdot \vol(B(R_i)) }
\end{equation*}
is monotone increasing in $i\ge 1$. By \eqref{e:lim vol>t}, the sequence $\varphi_{F}^{i}$ has the almost everywhere limit
\begin{equation}
\label{eq:lim phii=phi}
\quad \lim\limits_{i \rightarrow\infty} \varphi^i_F(t) = 1-\Psi_F(t) ,
\end{equation}
and by applying the Monotone Convergence Theorem on \eqref{eq:lim phii=phi}, we obtain the equality
\begin{equation}
\label{e:limits1}
\lim\limits_{i\rightarrow\infty} \int \limits_{0}^{\infty} \varphi^i_F(t)  \, dt = \int\limits_{0}^{\infty} (1 - \Psi_F(t)) \, dt.
\end{equation}
Observe that
\[  \int \limits_{0}^{\infty} \# \{ v\in V(R_i) : \vol(v) \ge t \}  \, dt  = \sum_{v \in V(R_i) }  \vol(v) ,\]
and so, interchanging expectation and integration,
\begin{equation}
\label{e:limits2}
\int \limits_{0}^{\infty} \varphi^i_F(t)  \, dt = \frac{ \E\left[ \sum_{v \in V(R_i) }  \vol(v) \right]  }{ c_{NS}(\rho)\cdot\vol(B(R_i)) }  .
\end{equation}
Combining \eqref{e:limits1} and \eqref{e:limits2}, we conclude that
\begin{equation}
\label{e:finallimits}
\int\limits_{0}^{\infty} ( 1 -  \Psi_F(t)) \, dt =  \lim\limits_{i\rightarrow\infty}   \frac{ \E\left[  \sum_{v \in V(R_i) }  \vol(v)  \right]  }{ c_{NS}(\rho)\cdot \vol(B(R_i)) } .
\end{equation}

It remains to analyse the r.h.s.\ of equality \eqref{e:finallimits}. To this end we notice that, using the definition of the boundary volume $\Vc(R)$ in Definition~\ref{d:bouvol}, for every $R>0$
\begin{equation}
\label{eq:intern vol=tot-intsct}
\sum_{v \in V(R) }  \vol(v) = \vol(B(R)) - \Vc(R) .
\end{equation}
Inserting \eqref{eq:intern vol=tot-intsct} into \eqref{e:finallimits} finally yields
\begin{align*}
\int\limits_{0}^{\infty} (1 - \Psi_F(t)) \, dt & =     \lim\limits_{i\rightarrow\infty}   \frac{ \vol(B(R_i))  - \E[\Vc(R_i)]}{ c_{NS}(\rho)\cdot \vol(B(R_i)) } \\
&=     \frac{1}{c_{NS}(\rho)} \left( 1 - \lim\limits_{i\rightarrow\infty}   \frac{ \E[  \Vc(R_i)] }{ \vol(B(R_i)) }  \right) ,
\end{align*}
completing the proof of \eqref{eq:lim eq vol}.

We turn to statement \eqref{eq:lim eq nest}, which is proved similarly. Arguing as for the first statement, and replacing integrals with sums whenever necessary, we arrive at the following analogue of the equality \eqref{e:finallimits} for the connectivity:
 \begin{equation}
 \label{e:finallimits2}
  \sum\limits_{k=0}^{\infty} k \, \mu_{\Gamma(F)}(k)  =     \lim\limits_{i\rightarrow\infty}   \frac{ \E[  \sum_{v \in V(R_i)} d(v) ]  }{ c_{NS}(\rho)\vol(B(R_i)) }  .
  \end{equation}
Recall \eqref{e:tree}, which states that
\[ \sum_{v \in \widebar{V}(R)} \widebar{d}(v) = 2(\widebar{\Nc}(F; R) - 1)  .\]
By the definition of the boundary connectivity $\Ccc(R)$ in Definition~\ref{def:bndr conn}, we therefore have
\begin{align}
\label{eq:intern vol=tot-intsct2}
 \sum_{v \in V(R)}d(v) &=  2(\widebar{\Nc}(F; R) - 1) - \Ccc(R) \\
 \nonumber &  = 2(\Nc(F; R) - 1) - \Ccc(R) + 2(\widebar{\Nc}(F; R) - \Nc(F; R))  .
 \end{align}
Inserting \eqref{eq:intern vol=tot-intsct2} into \eqref{e:finallimits2} yields
\[  \sum\limits_{k=0}^{\infty} k \,\mu_{\Gamma(F)}(k)  =     \lim\limits_{i\rightarrow\infty}   \frac{ 2\E[ \Nc(F; R_i) ] - \E[\Ccc(R_i)] + 2 \E[\widebar{\Nc}(F; R) - \Nc(F; R)] }{ c_{NS}(\rho)\vol(B(R_i)) }  .   \]
Given Proposition \ref{p:boun} and the convergence in \eqref{eq:NS conv mean}, this reduces to
\[  \sum\limits_{k=0}^{\infty} k \, \mu_{\Gamma(F)}(k)  =    2 -  \frac{1}{c_{NS}(\rho)} \cdot \lim\limits_{i\rightarrow\infty}   \frac{ \E[  \Ccc(R_i)] }{ \vol(B(R_i)) }  , \]
completing the proof of \eqref{eq:lim eq nest}.
\end{proof}

\bigskip
\section{The empirical mean volume}
\label{s:emp}

Recall that Nazarov--Sodin showed \eqref{eq:NS conv mean} that, under the assumptions $(\rho 1)$--$(\rho 4)$,
\[  \frac{\Nc(F;R)}{\vol B(R)}  \to c_{NS}(\rho) \quad \text{in mean.} \]
In this section we verify that, under the additional nodal lower concentration property \eqref{eq:lower conc dec}, the `reciprocal' convergence
\[   \frac{\vol B(R)}{\Nc(F;R)} \id_{ \Nc(F;R) > 0 } \to \frac{1}{c_{NS}(\rho)} \quad \text{in mean}  \]
holds, i.e.\ the `empirical volume mean' converges to $1/c_{NS}(\rho)$ (see Theorem \ref{thm:mean cons vol Eucl}\eqref{it:empirical mean=1/NS}). The proof of Theorem \ref{thm:mean cons vol Eucl}\eqref{it:empirical mean=1/NS} only uses elementary properties of convergence in mean, and the proof of the related Theorem \ref{thm:mean cons vol Riem} is similar.

\begin{proof}[Proof of Theorem \ref{thm:mean cons vol Eucl}\eqref{it:empirical mean=1/NS}]
For every fixed $\varepsilon \in (0,1/c_{NS}(\rho))$ we may write
\begin{equation}
\label{eq:E=E1+...+E4}
\E \left[ \left|  \frac{\vol B(R)}{\Nc(F;R)}  \id_{\Nc(F;R) > 0 } - \frac{1}{c_{NS}(\rho)} \right| \right] =E_{1}+E_{2}+E_{3}+E_{4} ,
\end{equation}
where
\[ E_1 = E_{1}(F;R)=\E \left[ \left|    \frac{\vol B(R)}{\Nc(F;R)}  - \frac{1}{c_{NS}(\rho)} \right| \id_{ \vol(B(R)/\Nc(F;R) < 1/c_{NS}(\rho) - \varepsilon} \right]  , \]
\[ E_2 = E_{2}(F;R)=\E \left[ \left|    \frac{\vol B(R)}{\Nc(F;R)}  - \frac{1}{c_{NS}(\rho)} \right| \id_{ \vol(B(R)/\Nc(F;R) \in [ 1/c_{NS}(\rho) - \varepsilon, 1/c_{NS}(\rho) + \varepsilon ]} \right] ,\]
\[E_3 = E_{3}(F;R)=\E \left[ \left|    \frac{\vol B(R)}{\Nc(F;R)}  - \frac{1}{c_{NS}(\rho)} \right| \id_{ 1/c_{NS}(\rho) - \varepsilon <  \vol(B(R))/\Nc(F;R) < \infty} \right]  \]
and
\[E_4 = E_{4}(F;R)=\frac{1}{c_{NS}(\rho)} \prob( \Nc(F;R) = 0 ) .\]
Next we bound each of the $E_{i}$, $i=1,\ldots,4$ separately.

First,
\begin{equation}
\label{eq:E1->0}
E_{1}\le \frac{1}{c_{NS}(\rho)} \times \prob \left( \frac{\Nc(F;R)}{\vol(B(R)} > 1 / (1/c_{NS}(\rho) - \varepsilon ) \right) \rightarrow 0
\end{equation}
as $R \to \infty$, by \eqref{eq:NS conv mean}, (see, e.g., the law of large numbers \eqref{eq:NS law large numb}). Second, trivially
\begin{equation}
\label{eq:E2->0}
E_{2} <\epsilon.
\end{equation}
Next, since, being an integer, $\Nc(F;R)\ge 1$, we have
\begin{equation}
\label{eq:E3->0}
E_{3}\le  \max\{ \vol B(R), 1/c_{NS} \} \times \prob \left( \frac{\Nc(F;R)}{\vol(B(R))} < 1 / (1/c_{NS}(\rho) - \varepsilon ) \right)
\rightarrow 0
\end{equation}
by the definition \eqref{eq:lower conc dec} of nodal lower concentration, and $\vol B(R) =O(R^{d})$.
Lastly,
\begin{equation}
\label{eq:E4->0}
E_4 \rightarrow 0,
\end{equation}
since
\[  \prob( \Nc(F;R) = 0 ) \to 0  \]
by the law of large numbers \eqref{eq:NS law large numb}.
We finally collect \eqref{eq:E1->0}, \eqref{eq:E2->0}, \eqref{eq:E3->0} and \eqref{eq:E4->0}, substitute these into \eqref{eq:E=E1+...+E4},
and take $\varepsilon \to 0$, to establish that
\[  \E \left[ \left|    \frac{\vol B(R)}{\Nc(F;R)}  \id_{ \Nc(F;R) > 0 } - \frac{1}{c_{NS}(\rho)} \right|  \right]  \to 0 \]
as $R \to \infty$, which is the statement of Theorem \ref{thm:mean cons vol Eucl}\eqref{it:empirical mean=1/NS}.
\end{proof}

The proof of Theorem \ref{thm:mean cons vol Riem} is almost identical to the above:

\begin{proof}[Proof of Theorem \ref{thm:mean cons vol Riem}]
The statement \eqref{eq:empirical mean=1/NS ens} of Theorem \ref{thm:mean cons vol Riem} follows from the same argument as presented within the proof of Theorem \ref{thm:mean cons vol Eucl}\eqref{it:empirical mean=1/NS} above, where we replace \eqref{eq:NS conv mean} with its manifold version \eqref{eq:NS conv mean man}, and the nodal lower concentration property in Definition \eqref{def:nodal low conc} with its manifold version in Definition \ref{def:nod low conc ens}.
\end{proof}


\medskip
\bibliographystyle{alpha}
\bibliography{paper}

\newcommand{\etalchar}[1]{$^{#1}$}
\begin{thebibliography}{PCR{\etalchar{+}}98}

\bibitem[Ale96]{Al}
K.S. Alexander.
\newblock Boundedness of level lines for two-dimensional random fields.
\newblock {\em Ann. Probab.}, 24:1653--1674, 1996.

\bibitem[AW09]{AzWsh}
J.-M. Aza{\"{\i}}s and M.~Wschebor.
\newblock {\em Level sets and extrema of random processes and fields}.
\newblock John Wiley \& Sons, Inc., Hoboken, NJ, 2009.

\bibitem[BE87]{Bond}
J.R. Bond and G.~Efstathiou.
\newblock The statistics of cosmic background radiation fluctuations.
\newblock {\em Mon. Notices Royal Astron. Soc.}, 226(3):655--687, 1987.

\bibitem[Ber77]{Berry77}
M.V. Berry.
\newblock Regular and irregular semiclassical wavefunctions.
\newblock {\em J. Phys. A.}, 10(12):2083, 1977.

\bibitem[BG17]{BeGa}
V.~Beffara and D.~Gayet.
\newblock Percolation of random nodal lines.
\newblock {\em Publ. Math. IHES}, 126:131--176, 2017.

\bibitem[BLM87]{BrLeMa}
J.~Bricmont, J.L. Lebowitz, and C.~Maes.
\newblock Percolation in strongly correlated systems: the massless {G}aussian
  field.
\newblock {\em Jour. Stat. Phys.}, 48(5--6):1249--1268, 1987.

\bibitem[BMR18]{BeMuRi}
D.~Beliaev, S.~Muirhead, and A.~Rivera.
\newblock A covariance formula for topological events of smooth {G}aussian
  fields.
\newblock {\em arXiv preprint arxiv:1811.08169}, 2018.

\bibitem[BMW17]{BMW}
D.~Beliaev, S.~Muirhead, and I.~Wigman.
\newblock Russo-{S}eymour-{W}elsh estimates for the {K}ostlan ensemble of
  random polynomials.
\newblock {\em arXiv preprint arXiv:1709.08961}, 2017.

\bibitem[BR06]{BoRi}
B.~Bollob{\'a}s and O.~Riordan.
\newblock {\em Percolation}.
\newblock Cambridge University Press, 2006.

\bibitem[BS02]{BS}
E.~Bogomolny and C.~Schmit.
\newblock Percolation model for nodal domains of chaotic wave functions.
\newblock {\em Phys. Rev. Lett.}, 88(11):114102, Mar 2002.

\bibitem[BW18]{BeWi}
D.~Beliaev and I.~Wigman.
\newblock Volume distribution of nodal domains of random band-limited
  functions.
\newblock {\em Probab. Theory Relat. Fields}, 172(1--2):453--492, 2018.

\bibitem[Cil93]{Ci}
J.~Cilleruelo.
\newblock The distribution of the lattice points on circles.
\newblock {\em J. Number Theory}, 43(2):198202, 1993.

\bibitem[CS14]{CS}
Y.~Canzani and P.~Sarnak.
\newblock On the topology of the zero sets of monochromatic random waves.
\newblock {\em arXiv preprint arXiv:1412.4437}, 2014.

\bibitem[DPR18]{DrPrRo}
A.~Drewitz, A.~Pr{\'{e}}vost, and P.F. Rodriguez.
\newblock The sign clusters of the massless {G}aussian free field percolate on
  ${\mathbb{z}^{d}, d\geqslant 3}$ (and more).
\newblock {\em Commun. Math. Phys.}, 362(1), 2018.

\bibitem[GW11]{GaWe}
D.~Gayet and J.-Y. Welschinger.
\newblock Exponential rarefaction of real curves with many components.
\newblock {\em Publ. Math. IHES}, 113(1):69--96, 2011.

\bibitem[H\"68]{Ho}
L.~H\"{o}rmander.
\newblock The spectral function of an elliptic operator.
\newblock {\em Acta Math.}, 121:193--218, 1968.

\bibitem[KKW13]{KKW}
M.~Krishnapur, P.~Kurlberg, and I.~Wigman.
\newblock Nodal length fluctuations for arithmetic random waves.
\newblock {\em Ann. Math.}, 177(2):699--737, 2013.

\bibitem[KW17]{KuWi1}
P.~Kurlberg and I.~Wigman.
\newblock On probability measures arising from lattice points on circles.
\newblock {\em Math. Ann.}, 367(3--4):1057--1098, 2017.

\bibitem[KW18]{KuWi2}
P.~Kurlberg and I.~Wigman.
\newblock Variation of the nazarov-sodin constant for random plane waves and
  arithmetic random waves.
\newblock {\em Adv. Math.}, 330:516--552, 2018.

\bibitem[KZ14]{KlZi}
P.~Kleban and R.~Ziff.
\newblock Notes on connections in percolation clusters.
\newblock Private communication, 2014.

\bibitem[Lax57]{La}
P.D. Lax.
\newblock Asymptotic solutions of oscillatory initial value problems.
\newblock {\em Duke Math. J.}, 24:627--646, 1957.

\bibitem[LH57]{Longuet}
M.S. Longuet-Higgins.
\newblock The statistical analysis of a random, moving surface.
\newblock {\em Phil. Trans. R. Soc. Lond. A}, 249(966):321--387, 1957.

\bibitem[LPS09]{LPS}
H.~Lapointe, I.~Polterovich, and Y.~Safarov.
\newblock Average growth of the spectral function on a {R}iemannian manifold.
\newblock {\em Commun. Part. Diff. Eq.}, 34:581--615, 2009.

\bibitem[MS83a]{MoSt1}
S.A. Molchanov and A.K. Stepanov.
\newblock Percolation in random fields. {I}.
\newblock {\em Theor. Math. Phys.}, 55(2):478--484, 1983.

\bibitem[MS83b]{MoSt2}
S.A. Molchanov and A.K. Stepanov.
\newblock Percolation in random fields. {II}.
\newblock {\em Theor. Math. Phys.}, 55(3):592--599, 1983.

\bibitem[MV18]{MuVa}
S.~Muirhead and H.~Vanneuville.
\newblock The sharp phase transition for level set percolation of smooth planar
  gaussian fields.
\newblock {\em arXiv preprint arXiv:1806.11545}, 2018.

\bibitem[NS09]{NaSo09}
F.~Nazarov and M.~Sodin.
\newblock On the number of nodal domains of random spherical harmonics.
\newblock {\em Amer. J. Math.}, 131(5):1337--1357, 2009.

\bibitem[NS16]{NaSo15}
F.~Nazarov and M.~Sodin.
\newblock Asymptotic laws for the spatial distribution and the number of
  connected components of zero sets of {G}aussian random functions.
\newblock {\em J. Math. Phys. Anal. Geo.}, 12(3):205--278, 2016.

\bibitem[PCR{\etalchar{+}}98]{Park}
C.~Park, W.N. Colley, B.~Ratra, N.~Spergel, and N.~Sugiyama.
\newblock Cosmic microwave background anisotropy correlation function and
  topology from simulated maps for {MAP}.
\newblock {\em Astrophys. J.}, 506(2):473, 1998.

\bibitem[Ric44]{Rice}
S.O. Rice.
\newblock Mathematical analysis of random noise.
\newblock {\em Bell System Technical Journal}, 23(3):282--332, 1944.

\bibitem[Roz16]{Ro}
Y.~Rozenshein.
\newblock The number of nodal components of arithmetic random waves.
\newblock {\em Int. Math. Res. Notices.}, 2016.

\bibitem[RV17]{RiVa}
A.~Rivera and H.~Vanneuville.
\newblock Quasi-independence for nodal lines.
\newblock {\em arXiv preprint arXiv:1711.05009}, 2017.

\bibitem[Sar17]{Sa}
P.~Sarnak.
\newblock Private communication, 2017.

\bibitem[Sod12]{So12}
M.~Sodin.
\newblock Lectures on random nodal portraits.
\newblock {\em Lecture notes for a mini-course given at the St. Petersburg
  Summer School in Probability and Statistical Physics (June, 2012) Available
  at: http://www.math.tau.ac.il/sodin/SPB-Lecture-Notes.pdf}, 2012.

\bibitem[SS93]{ShSm}
M.~Shub and S.~Smale.
\newblock Complexity of {B}ezout's theorem. {II}. {V}olumes and probabilities.
\newblock In {\em Computational Algebraic Geometry. Progress in Mathematics,
  vol 109.} Birkhäuser, Boston, 1993.

\bibitem[SW15]{SaWi}
P.~Sarnak and I.~Wigman.
\newblock Topologies of nodal sets of random band limited functions.
\newblock {\em arXiv preprint arXiv:1510.08500}, 2015.

\bibitem[Swe62]{Swerling}
P.~Swerling.
\newblock Statistical properties of the contours of random surfaces.
\newblock {\em IEEE Trans. Inform. Theory}, 8(4):315--321, 1962.

\bibitem[Szn10]{Sz1}
A.S. Sznitman.
\newblock Vacant set of random interlacements and percolation.
\newblock {\em Ann. Math.}, 171(3):2039--2087, 2010.

\end{thebibliography}

\end{document}